\documentclass[12pt]{article}
\usepackage{a4wide,epsfig}
\usepackage{multirow}
\usepackage{amsfonts,amsmath}
\usepackage{subcaption}
\usepackage{booktabs}
\usepackage{pgfplotstable}

\newtheorem{theorem}{Theorem}[section]

\newtheorem{corollary}[theorem]{Corollary}
\newtheorem{lemma}[theorem]{Lemma}

\newtheorem{remark}{Remark}[section]
\newcommand{\proof} [1]
   { \noindent {\bf Proof.} #1 \hfill\rule{0.5em}{1.2ex} \par\medskip}

\numberwithin{equation}{section} 

\begin{document}

\setcounter{page}{1}

\title{Space-time finite element methods for
  distributed optimal control of the wave equation}
\author{Richard~L\"oscher, Olaf~Steinbach}
\date{Institut f\"ur Angewandte Mathematik, TU Graz, \\[1mm]
Steyrergasse 30, 8010 Graz, Austria}

\maketitle

\begin{abstract}
  We consider space-time tracking type distributed optimal control
  problems for the wave equation in the space-time domain
  $Q:= \Omega \times (0,T) \subset {\mathbb{R}}^{n+1}$,
  where the control is assumed to be in the energy space
  $[H_{0;,0}^{1,1}(Q)]^*$, rather than in $L^2(Q)$ which is more common.
  While the latter ensures a unique state in the Sobolev space
  $H^{1,1}_{0;0,}(Q)$, this does not define a solution isomorphism.
  Hence we use an appropriate state space $X$ such that
  the wave operator becomes an isomorphism from $X$ onto
  $[H_{0;,0}^{1,1}(Q)]^*$. Using space-time finite element spaces of
  piecewise linear continuous basis functions on completely unstructured
  but shape regular simplicial meshes, we derive a priori estimates for
  the error $\|\widetilde{u}_{\varrho h}-\overline{u}\|_{L^2(Q)}$ between the
  computed space-time finite element solution $\widetilde{u}_{\varrho h}$ and
  the target function $\overline{u}$ with respect to the regularization
  parameter $\varrho$, and the space-time finite element mesh-size $h$,
  depending on the regularity of the desired state $\overline{u}$.
  These estimates lead to the optimal choice $\varrho=h^2$ in order to
  define the regularization parameter $\varrho$ for a given space-time
  finite element mesh size $h$, or to determine the required mesh size $h$
  when $\varrho$ is a given constant representing the costs of the control.
  The theoretical results will be supported by numerical examples with
  targets of different regularities, including discontinuous targets.
  Furthermore, an adaptive space-time finite element scheme is proposed
  and numerically analyzed.                                      
\end{abstract}


\section{Introduction}
We consider a distributed optimal control
problem to minimize a tracking type functional to reach a given
target $\overline{u} \in L^2(Q)$ subject to the
initial boundary value problem for the wave equation with zero
initial and boundary conditions in the space-time domain $Q$. The standard
setting of such kind of optimal control problems assumes the control
to be in $L^2(Q)$, see, e.g.,
\cite{LionsVol1:1971,KunischPeralta:2022,Troeltzsch:2010}.
In this case, the wave equation admits a unique solution in the Sobolev space
$H_{0;0,}^{1,1}(Q)$, see \cite{Ladyzhenskaya:1985,SteinbachZank:2020}.
For our analysis though, we will use a regularization in the (energy)
space $[H_{0;,0}^{1,1}(Q)]^\ast$ which is the dual of the test space
for the variational formulation of the wave equation. To ensure unique
solvability of the wave equation also in this case, we use a generalized
variational formulation of the wave equation as recently discussed in
\cite{SteinbachZank:2022}. Similar investigations using the energy norm
for the control were already done for distributed optimal control problems
subject to elliptic
\cite{LLSY2:LangerSteinbachYang:2022b,LLSY2:NeumuellerSteinbach:2021a}
and parabolic partial differential equations
\cite{LangerSteinbachTroeltzschYang:2021,UL:LangerSteinbachYang:2022a},
and, as it turns out, our analysis fits into the same framework. 

In this paper, our main interest will be in proving estimates for the error
$\|\widetilde{u}_{\varrho h}-\overline u\|_{L^2(Q)}$ for the computed
space-time finite element solution $\widetilde{u}_{\varrho h}$,
depending on the regularity of the target function $\overline{u}$
and on the regularization parameter $\varrho$. In particular, in the
discrete setting, we will allow $\varrho$ to depend on the mesh size $h$
and we derive an optimal choice $\varrho=h^2$ in the sense, that we can
achieve optimal orders of convergence with respect to the regularity of
$\overline{u}$. This is of particular interest when the regularization
parameter $\varrho$ is required to ensure solvability of the unconstrained
optimal control problem, i.e., the costs are not of practical interest,
see, e.g., \cite{Poerner:2018}. 
In this case, the minimization problem is
closely related to the Tikhonov regularization in inverse problems,
where the parameter dependent convergence as $\varrho \to 0$ is
well-studied, see, e.g.,
\cite{AlbaniCezaroZubelli:2016,EnglHankeNeubauer:1996,Isakov:2017}. 
On the other hand, when $\varrho$ is a given constant representing
the costs of the control, one can determine the required space-time finite
element mesh size $h$ in order to reach the minimum of the functional
to be minimized. The optimal relation between the regularization
parameter $\varrho$ and the finite element mesh size $h$ is also important
for the design of preconditioned iterative solution methods for the
discrete optimality system, see, e.g.,
\cite{UL:LangerSteinbachYang:2022a,LLSY2:LangerSteinbachYang:2022b}
for the elliptic and the parabolic case, respectively.
To ease the presentation, at this time, we will not
consider any control or state constraints, see, e.g.,
\cite{GugatKeimerLeugering:2009,KroenerKunischVexler:2011}.
However, state or control constraints can be considered within the
abstract framework as given in \cite{GanglLoescherSteinbach:2022}.

When choosing an appropriate state space $X$ as introduced in
\cite{SteinbachZank:2022}, the state equation, i.e., the Dirichlet
problem for the wave equation, admits a unique solution $u_\varrho \in X$,
for each right hand side $z_\varrho \in Z = Y^* = [H_{0;,0}^{1,1}(Q)]^*$,
i.e., the wave operator $ B: X \to Y^*$ is an isomorphism.
In view of the Ne\u{c}as--Babu\v{s}ka theorem, e.g.,
\cite{LSTY:BabuskaAziz:1972a,LSTY:Necas:1962a},
$B$ is, in particular, inf-sup stable. Furthermore, when introducing
a self-adjoint, elliptic and bounded operator $A:Y\to Y^*$, which gives
raise to an equivalent norm in $Y$, and  $p_\varrho \in Y$ as the
solution of the adjoint wave equation $B^* p_\varrho= u_\varrho-\overline u$,
we can eliminate the control $z_\varrho$ by the gradient equation
$p_\varrho+\varrho A^{-1}z_\varrho=0$. Then, the unique solution of the
optimal control problem can be computed by solving the reduced first
order optimality system
\[
\begin{pmatrix}
	\varrho^{-1} A & B\\
	-B^* & I 
\end{pmatrix}
\begin{pmatrix}
	p_\varrho\\
	u_\varrho
\end{pmatrix}
=
\begin{pmatrix}
	0\\
	\overline u
\end{pmatrix},
\]
for any given target function $\overline{u} \in L^2(Q)$, which can be
interpreted as a stabilized saddle point formulation. This specific form
arises also in boundary optimal control problems for the wave equation,
see, e.g., \cite{MontanerMunch:2019}, and, undoubtely, in many applications.

For the numerical treatment of the above considered optimal control problem,
there are a myriad of methods available, e.g.,
\cite{GlowinskiKintonWheeler:1989,HulbertHughes:1990,
	KroenerKunischVexler:2011,KunischPeralta:2022,Zuazua:2005},
just to mention a few. Mostly space and time are treated separately, using,
e.g., finite difference methods, mixed and discontinuous Galerkin finite
element methods, finite volume methods, and time stepping schemes or
variational in time methods. Here, we will consider a real space-time
finite element method on completely unstructured, but shape regular,
simplicial space-time finite element meshes decomposing the space-time
domain $Q$, see also
\cite{DoerflerFindeisenWienersZiegler:2019,ErnestiWieners:2019} where such
methods are given for the direct solution of the wave equation.
Introducing conforming finite element spaces $X_h \subset X$ and
$Y_h\subset Y$ with appropriate approximation properties, the discrete
reduced optimality system admits again the form of a stabilized saddle
point formulation. Though, at this point it is worth stressing, that the
assumptions on the discrete operator $B_h : X_h\to Y_h^*$ are vastly
weakened, i.e., we do not need a discrete inf-sup stability condition and
not even a CFL-condition to be fulfilled. Therefore, this method directly
allows for an adaptive finite element scheme, see, e.g.,
\cite{MeidnerVexler:2007} in the case of a parabolic optimal control problem,
and \cite{BangerthGeigerRannacher:2010} for adaptive schemes for the wave
equation, which we will also address in our numerical investigations. 

The remainder of this paper is structured as follows: In Section 2 we will
state the model problem and introduce the appropriate functional analytical
setting required for the solution of the wave equation. In Section 3 we present
the main result for the regularization error estimates which depend on the
regularity of the target $\overline u$, and on the regularization
parameter $\varrho$. The space-time finite element discretization and
related a priori error estimates are investigated in Section 4, where we
will conclude the optimal choice $\varrho=h^2$ for the regularization parameter.
Numerical tests will confirm our theory in Section 5. Furthermore, we will
compare the proposed energy regularization approach with the more
standard $L^2$ regularization in the same setting, as well as propose
an adaptive refinement strategy. In Section 6, we draw some conclusions
and give an outlook on ongoing work. 

\section{Distributed optimal control problems}
Let $\Omega \subset {\mathbb{R}}^n$, $n=1,2,3$, be a bounded convex domain with,
for $n=2,3$, Lipschitz boundary $\Gamma = \partial \Omega$, and let
$T >0$ be a given finite time horizon. Then we introduce the space-time
domain $Q := \Omega \times (0,T)$ and the lateral boundary
$\Sigma := \Gamma \times (0,T)$.
For a given target $\overline{u} \in L^2(Q)$ and a regularization parameter
$\varrho > 0$, we consider the minimization of the cost functional
\begin{equation}\label{cost functional}
	{\mathcal{J}}(u_\varrho,z_\varrho) := \frac{1}{2} \int_0^T \int_\Omega
	[u_\varrho(x,t)-\overline{u}(x,t)]^2 \, dx \, dt
	+
	\frac{1}{2} \, \varrho \, \| z_\varrho \|^2_Z
\end{equation}
subject to the initial boundary value problem for the wave equation with
homogeneous Dirichlet boundary conditions,
\begin{equation}\label{DBVP wave equation}
	\begin{array}{rclcl}
		\Box u_\varrho(x,t) :=
		\partial_{tt} u_\varrho(x,t) - \Delta_x u_\varrho(x,t)
		& = & z_\varrho(x,t)
		& & \mbox{for} \; (x,t) \in Q, \\[1mm]
		u_\varrho(x,t)
		& = & 0
		& & \mbox{for} \; (x,t) \in \Sigma,
		\\[1mm]
		u_\varrho(x,0) = \partial_t u_\varrho(x,t)_{|t=0}
		& = & 0 & & \mbox{for} \; x \in \Omega .
	\end{array}
\end{equation}
Our particular interest is in the numerical solution of the
constrained minimization problem \eqref{cost functional} and
\eqref{DBVP wave equation} by using a space-time finite element approach
on simplicial meshes which are completely unstructured in space and time.
For the error $\| \widetilde{u}_{\varrho h} - \overline{u} \|_{L^2(Q)}$ of the
computed numerical solution $\widetilde{u}_{\varrho h}$ we will provide
estimates in the
space-time finite element mesh size $h$, and in the regularization
parameter $\varrho$ from which we will
derive an optimal choice for $\varrho$, which will depend on the
choice of the regularization space $Z$.

First we consider $z_\varrho \in Z = L^2(Q)$. Following \cite{SteinbachZank:2020},
the space-time variational formulation of the state equation
\eqref{DBVP wave equation}
is to find $ u_\varrho \in H^{1,1}_{0;0,}(Q) $ such that
\begin{equation}\label{DBVP wave equation VF}
	b(u_\varrho,q) :=
	- \langle \partial_t u_\varrho , \partial_t q \rangle_{L^2(Q)} +
	\langle \nabla_x u_\varrho , \nabla_x q \rangle_{L^2(Q)} =
	\langle z_\varrho , q \rangle_{L^2(Q)}
\end{equation}
is satisfied for all $q \in H^{1,1}_{0;,0}(Q)$. Here we use the
anisotropic Sobolev space
\[
H^{1,1}_{0;0,}(Q) := L^2(0,T;H^1_0(\Omega)) \cap H^1_{0,}(0,T;L^2(\Omega)),
\]
where $H^1_{0,}(0,T;L^2(\Omega))$ covers the zero initial condition
$u(x,0)=0$ for $x \in \Omega$, while $L^2(0,T;H^1_0(\Omega))$ includes
the homogeneous Dirichlet boundary condition on $\Sigma$. Note that the
second initial condition $\partial_t u(x,t)_{|t=0}=0$ for $x \in \Omega$
enters the variational formulation \eqref{DBVP wave equation VF} in
a natural way. A norm in $H^{1,1}_{0;0,}(Q)$ is given by the graph norm
\[
\| u \|_{H^{1,1}_{0;0,}(Q)} := \sqrt{ \| \partial_t u \|_{L^2(Q)}^2 +
	\| \nabla_x u \|^2_{L^2(Q)}} \, = \, |u|_{H^1(Q)} .
\]
Note that $H^{1,1}_{0;,0}(Q)$ is defined accordingly, but with a zero
terminal condition $q(x,T)=0$ for $x \in \Omega$. Then we have
\begin{equation}\label{Boundedness b}
	|b(u,q)| \, \leq \, |u|_{H^1(Q)} |q|_{H^1(Q)} \quad
	\mbox{for all} \; u \in H^{1,1}_{0;0,}(Q), q \in H^{1,1}_{0;,0}(Q) .
\end{equation}
For $z_\varrho \in L^2(Q)$ there exists a unique solution
$u_\varrho \in H^{1,1}_{0;0,}(Q)$ of the variational formulation
\eqref{DBVP wave equation VF} satisfying,
see, e.g., \cite[Theorem 5.1]{SteinbachZank:2020}, and
\cite{Ladyzhenskaya:1985},
\[
\| u_\varrho \|_{H^{1,1}_{0;0,}(Q)} \leq \frac{1}{\sqrt{2}} \, T \,
\| z_\varrho \|_{L^2(Q)} \, .
\]
Hence we can write $u_\varrho = {\mathcal{S}} z_\varrho$ with the
solution operator ${\mathcal{S}} : L^2(Q) \to H^{1,1}_{0;0,}(Q) \subset L^2(Q)$,
and we can introduce the reduced cost functional
\[
\widetilde{\mathcal{J}}(z_\varrho) :=
\frac{1}{2} \, \| {\mathcal{S}} z_\varrho - \overline{u} \|^2_{L^2(Q)} +
\frac{1}{2} \, \varrho \, \| z \|^2_{L^2(Q)},
\]
whose minimizer is given by the gradient equation
\begin{equation}\label{Gradient equation}
	p_\varrho(x,t) + \varrho \, z_\varrho(x,t) = 0 \quad \mbox{for} \; (x,t) \in Q,
\end{equation}
and where $p_\varrho \in H^{1,1}_{0;,0}(Q)$ is the weak solution of the adjoint
problem
\begin{equation}\label{DBVP adjoint wave equation}
	\begin{array}{rclcl}
		\partial_{tt} p_\varrho(x,t) - \Delta_x p_\varrho(x,t)
		& = & u_\varrho(x,t) - \overline{u}(x,t)
		& & \mbox{for} \; (x,t) \in Q, \\[1mm]
		p_\varrho(x,t)
		& = & 0
		& & \mbox{for} \; (x,t) \in \Sigma,
		\\[1mm]
		p_\varrho(x,T) = \partial_t p_\varrho(x,t)_{|t=T}
		& = & 0 & & \mbox{for} \; x \in \Omega .
	\end{array}
\end{equation}
Similar as in \cite{LangerSteinbachTroeltzschYang:2021} for the heat equation we can apply a
space-time finite element method on completely unstructured simplicial
meshes to discretize the optimality system \eqref{DBVP wave equation}
and \eqref{DBVP adjoint wave equation} after eliminating the control
$z_\varrho$ from the gradient equation \eqref{Gradient equation}.
Although we will consider this approach for a numerical comparison,
at this time we are not able to provide a complete numerical analysis
for this approach. As already seen in the elliptic case
\cite{LLSY2:LangerSteinbachYang:2022b, LLSY2:NeumuellerSteinbach:2021a},
and in the parabolic case
\cite{LangerSteinbachTroeltzschYang:2021,LangerSteinbachTroeltzschYang:2022,
	UL:LangerSteinbachYang:2022a}, there are differences
both in the numerical analysis and in the properties of the numerical
solutions when considering the regularization in $L^2(Q)$, and in the
related energy space, which is the dual of the test space.

A direct space-time finite element discretization
of the variational formulation \eqref{DBVP wave equation VF} on space-time
tensor product meshes using piecewise linear continuous basis functions
requires an appropriate stability condition
$h_t \leq h_x/\sqrt{n}$ where $h_t$ and $h_x$ are the temporal
and spatial mesh sizes, respectively, see \cite{SteinbachZank:2020}.
Moreover, the associated operator $B$ to the variational formulation
\eqref{DBVP wave equation VF} does not define
an isomorphism between $L^2(Q)$ and $H^{1,1}_{0;0,}(Q)$, see
Theorem \ref{Theorem 2.1}.
Although the variational formulation \eqref{DBVP wave equation VF} is well
defined also for $z_\varrho \in [H^{1,1}_{0;,0}(Q)]^*$, it does not ensure
unique solvability in $H^{1,1}_{0;0,}(Q)$ in this case. Instead we have to
enlarge the ansatz space in order to
incorporate the second initial condition
$\partial_t u_\varrho(x,t)_{|t=0}=0$ in an appropriate way.
In what follows we will consider a generalized variational formulation
of the wave equation, see \cite{SteinbachZank:2022}. When
using a distributional definition of the wave operator, we consider
an ultra-weak variational formulation of \eqref{DBVP wave equation}
to find $u \in L^2(Q)$ which is extended by zero to an enlarged domain to
cover the initial conditions.
This approach will allow us
to define the regularization in a suitable energy norm. In this case we
choose $Z=[H_{0;,0}^{1,1}(Q)]^*$ as the dual of the test space. A norm in
this space is given as
\[
\| z \|_{[H_{0;,0}^{1,1}(Q)]^*} :=
\sup_{0\neq q\in H^{1,1}_{0;,0}(Q)}
\frac{\langle z ,q\rangle _Q}{\|q\|_{H^{1,1}_{0;,0}(Q)}},
\]
where $\langle \cdot, \cdot\rangle_Q$ is an extension of the inner product
in $L^2(Q)$. For $z \in [H_{0;,0}^{1,1}(Q)]^*$, and using the Riesz
isomorphism, there exists a unique $w_z \in H_{0;,0}^{1,1}(Q)$ such that 
\[
\langle A w_z , q \rangle_Q :=
\langle \partial_t w_z ,\partial_tq\rangle_{L^2(Q)} +
\langle \nabla_x w_z ,\nabla_x q\rangle_{L^2(Q)}
=
\langle z ,q\rangle_Q \quad \mbox{for all} \;
q \in H_{0;,0}^{1,1}(Q) .
\]
With this choice, $A$ is self-adjoint, elliptic and bounded, i.e., 
\[
\langle A w , q \rangle_Q \leq \, |w|_{H^1(Q)} | q |_{H^1(Q)}, \quad
\langle A q, q \rangle_Q = |q|^2_{H^1(Q)} \quad
\mbox{for all} \; w,q \in H^{1,1}_{0;,0}(Q) ,
\]
and, hence, invertible. Thus, we can write
\[
\| w_z \|^2_{H^{1,1}_{0;,0}(Q)} =
\langle A w_z , w_z \rangle_Q = \| z \|^2_{[H^{1,1}_{0;,0}(Q)]^*}
\]
as well as
\[
\|z \|_{[H_{0;,0}^{1,1}(Q)]^*}^2 =
\langle z , w_z \rangle_Q =
\langle z ,A^{-1}z \rangle _Q \quad \mbox{for all} \;
z \in [H^{1,1}_{0;,0}(Q)]^*.
\]
We proceed with stating some preliminaries. First, let us give a result
concerning the boundedness of the solution $u_\varrho \in H^{1,1}_{0;0,}(Q)$ of
\eqref{DBVP wave equation VF} when considering
the norm of $z_\varrho$ in $[H^{1,1}_{0;,0}(Q)]^*$.

\begin{theorem}{\rm \cite[Theorem 1.1]{SteinbachZank:2022}}
	\label{Theorem 2.1}
	There does not exist a constant $c>0$ such that each right-hand side
	$z_\varrho \in L^2(Q)$ and the corresponding solution
	$u_\varrho \in H^{1,1}_{0;0,}(Q)$ of \eqref{DBVP wave equation VF}
	satisfy
	\[
	\| u_\varrho \|_{H^{1,1}_{0;0,}(Q)} \leq c \,
	\|z_\varrho \|_{[H_{0;,0}^{1,1}(Q)]^*}.
	\]
	In particular, the inf-sup condition
	\[
	c_S \, \| u \|_{H^{1,1}_{0;0,}(Q)} \leq
	\sup_{0\neq q\in H^{1,1}_{0;,0}(Q)}
	\frac{b(u,q)}{\|q\|_{H^{1,1}_{0;,0}(Q)}} \quad
	\mbox{for all} \; u\in H^{1,1}_{0;0,}(Q)
	\]
	with a constant $c_S>0$ does not hold true. 
\end{theorem}  

\noindent
The issue to overcome is the handling of the initial condition
$\partial_t u_\varrho(x,t)_{|t=0}=0$ for $x \in \Omega$ for which we
will proceed as in \cite{SteinbachZank:2022}. For the
enlarged space-time domain $Q_-:=\Omega \times (-T,T)$, and for
$u\in L^2(Q)$ we define the zero extension
\[
\widetilde{u}(x,t) :=
\begin{cases}
	u(x,t) & \mbox{for} \; (x,t)\in Q, \\
	0, & \mbox{else.}
\end{cases}
\]
The application of the wave operator $\Box \widetilde{u}$ on $Q_-$ will
be formulated as a distribution, i.e., for
$\varphi \in C_0^\infty(Q_-)$ we define
\[
\langle \Box \widetilde{u}, \varphi \rangle_{Q_-}
:=
\int_{Q_-} \widetilde{u}(x,t) \, \Box \varphi(x,t) \, dx \, dt =
\int_Q u(x,t) \, \Box \varphi(x,t) \, dx \, dt.  
\]
Now we are in the position to introduce the space 
\[
\mathcal{H}(Q) := \Big \{
u=\widetilde{u}_{|_{Q}} : \widetilde{u} \in L^2(Q_-) , \,
\widetilde{u}_{|\Omega \times (-T,0)}=0, \,
\Box \widetilde{u} \in [H_0^1(Q_-)]^* \Big \},
\]
with the graph norm
\[
\|u\|_{\mathcal{H}(Q)} :=
\sqrt{\|u\|_{L^2(Q)}^2 + \|\Box \widetilde{u} \|_{[H_0^1(Q_-)]^*}^2} \, .
\]
The normed vector space $(\mathcal{H}(Q),\|\cdot\|_{\mathcal{H}(Q)})$ is a
Banach space, and it holds true that, see \cite[Lemma 3.5]{SteinbachZank:2022},
$H_{0;0,}^{1,1}(Q)\subset \mathcal{H}(Q)$ i.e.,
\begin{equation}\label{norm X H1}
	\| \Box \widetilde{u} \|_{[H^1_0(Q_-)]^*} \leq \| u \|_{H^{1,1}_{0;0,}(Q)}
	\quad \mbox{for all} \; u \in H^{1,1}_{0;0,}(Q).
\end{equation}
Therefore, we can consider the space
\[
\mathcal{H}_{0;0,}(Q) :=
\overline{H_{0;0,}^{1,1}(Q)}^{\| \cdot \|_{\mathcal{H}(Q)}} \subset
\mathcal{H}(Q)
\]
which will serve as ansatz space. For $u \in {\mathcal{H}}_{0;0,}(Q)$,
an equivalent norm is given as, see \cite[Lemma 3.6]{SteinbachZank:2022},
\[
\| u \|_{{\mathcal{H}}_{0;0,}(Q)} = \| \Box \widetilde{u}
\|_{[H^1_0(Q_-)]^*} .
\]
For given
$z_\varrho \in [H_{0;,0}^{1,1}(Q)]^*$ we consider the variational formulation
to find $u_\varrho \in \mathcal{H}_{0;0,}(Q)$ such that 
\begin{equation}\label{eq:generalized-VF-wave}
	\langle \Box \widetilde{u}_\varrho , {\mathcal{E}} q \rangle_{Q_-}
	= \langle z_\varrho , q \rangle_{Q} \quad
	\mbox{for all}  \; q\in H_{0;,0}^{1,1}(Q),
\end{equation}
where $\mathcal{E} : H_{0;,0}^{1,1}(Q) \to H_0^1(Q_-)$ is a suitable
extension operator, e.g., reflection in time with respect to $t=0$,
satisfying
\[
\| {\mathcal{E}} q \|_{H^1_0(Q_-)} \leq 2 \,
\| q \|_{H^{1,1}_{0;,0}(Q)} \quad \mbox{for all} \;
q \in H^{1,1}_{0;,0}(Q) .
\]
We conclude that the bilinear form within the variational
formulation \eqref{eq:generalized-VF-wave} is bounded, i.e., for
all $u \in {\mathcal{H}}_{0;0,}(Q)$ and $q \in H^{1,1}_{0;,0}(Q)$
we have
\begin{equation}\label{Boundedness Bilinear Form}
	\left| \langle \Box \widetilde{u}_\varrho , {\mathcal{E}} q \rangle_{Q_-}
	\right| \, \leq \,
	\| \Box \widetilde{u} \|_{[H^1_0(Q_-)]^*}
	\| {\mathcal{E}} q \|_{H^1_0(Q_-)} \leq 2 \,
	\| u \|_{{\mathcal{H}}_{0;0,}(Q)} \| q \|_{H^{1,1}_{0;,0}(Q)} .
\end{equation}
Moreover, we have the following result.

\begin{theorem}{\rm \cite[Theorem 3.9]{SteinbachZank:2022}}
	\label{thm:energy-reg_unique-solvability}
	For each given $z_\varrho \in [H_{0;,0}^{1,1}(Q)]^*$, there exists a unique
	solution $u_\varrho \in \mathcal{H}_{0;0,}(Q)$ of the variational
	formulation \eqref{eq:generalized-VF-wave} satisfying
	\[
	\| u_\varrho \|_{\mathcal{H}_{0;0,}(Q)} =
	\| \Box \widetilde{u}_\varrho\|_{[H_0^1(Q_-)]^*} =
	\| z_\varrho \|_{[H_{0;,0}^{1,1}(Q)]^*}.
	\]
	In particular, there holds the inf-sup stability condition
	\begin{equation}\label{inf-sup generalized}
		\| u \|_{\mathcal{H}_{0;0,}(Q)} \leq
		\sup_{0\neq q\in H_{0;,0}^{1,1}(Q)}
		\frac{\langle \Box \widetilde{u} , {\mathcal{E}} q \rangle_{Q_-}}
		{\| q \|_{H_{0;,0}^{1,1}(Q)}} \quad  \mbox{for all} \;
		u \in \mathcal{H}_{0;0,}(Q) .
	\end{equation}
\end{theorem} 

\begin{remark}\label{rem:regularity-bform_norm-estimate}
	The use of the bilinear form
	$\langle \Box \widetilde{u} , {\mathcal{E}} q \rangle_{Q_-}$
	might seem cumbersome. But we have, see
	{\rm \cite[Lemma 3.5]{SteinbachZank:2022}},
	\[
	\langle \Box \widetilde{u} , {\mathcal{E}} q \rangle_{Q_-}
	=
	- \langle \partial_t u , \partial_t q \rangle_{L^2(Q)} +
	\langle \nabla_x u , \nabla_x q \rangle_{L^2(Q)}
	\; \mbox{for all} \;
	u \in H^{1,1}_{0;0,}(Q) \subset {\mathcal{H}}_{0;0,}(Q),
	q \in H^{1,1}_{0;,0}(Q) .
	\]
	This is of particular interest when considering the discrete setting,
	as piecewise linear continuous functions are in $H^1(Q)$.
\end{remark}

\noindent
In view of Theorem \ref{thm:energy-reg_unique-solvability} we have a
solution operator $\mathcal{S} : [H_{0;,0}^{1,1}(Q)]^* \to
\mathcal{H}_{0;0,}(Q) \subset L^2(Q)$.
So, we can write the reduced cost functional
\begin{eqnarray*}
	\widetilde{\mathcal{J}}(z_\varrho)
	& := & \frac{1}{2} \, \| {\mathcal{S}} z_\varrho - \overline{u} \|_{L^2(Q)}^2
	+ \frac{1}{2} \, \varrho \, \| z_\varrho \|_{[H_{0;,0}^{1,1}(Q)]^*}^2 \\
	& = & \frac{1}{2} \, \langle {\mathcal{S}}^\ast {\mathcal{S}} z_\varrho,
	z_\varrho \rangle_Q - \langle {\mathcal{S}}^\ast \overline{u},
	z_\varrho \rangle _Q + \frac{1}{2} \, \| \overline{u} \|_{L^2(Q)}^2
	+ \frac{1}{2} \, \varrho \, \langle A^{-1} z_\varrho ,
	z_\varrho\rangle_Q,
\end{eqnarray*}
where $\mathcal{S}^\ast: [\mathcal{H}_{0;0,}(Q)]^* \to H_{0;,0}^{1,1}(Q)$
denotes the dual of the solution operator. The minimizer of the reduced
cost functional is the unique solution of the gradient equation
\begin{equation}\label{eq:gradient-energy}
	\mathcal{S}^\ast(\mathcal{S} z_\varrho - \overline u) +
	\varrho \, A^{-1} z_\varrho =0, 
\end{equation}
i.e., we have to find $z_\varrho \in [H_{0;,0}^{1,1}(Q)]^*$ as solution of
\[
\varrho \, A^{-1} z_\varrho + {\mathcal{S}}^\ast {\mathcal{S}} z_\varrho
= {\mathcal{S}}^\ast \overline{u} \quad \mbox{in} \;
H_{0;,0}^{1,1}(Q).
\]
Note that $T_\varrho := \varrho \, A^{-1}+ {\mathcal{S}}^\ast {\mathcal{S}} :
[H_{0;,0}^{1,1}(Q)]^* \to H_{0;,0}^{1,1}(Q)$ is bounded and elliptic, thus
unique solvability of the operator equation \eqref{eq:gradient-energy}
follows immediately. When introducing the adjoint state
$p_\varrho = {\mathcal{S}}^\ast(u_\varrho-\overline{u})$, and
$w_{z_\varrho} \in H_{0;,0}^{1,1}(Q)$ as solution of $Aw_{z_\varrho}=z_\varrho$,
we can write the gradient equation \eqref{eq:gradient-energy} as
\begin{equation}\label{gradient equation p z}
	p_\varrho + \varrho \, w_{z_\varrho} = 0,
\end{equation}
where $p_\varrho\in H_{0;,0}^{1,1}(Q)$ is the unique solution of the adjoint
generalized wave equation
\begin{equation}\label{eq:adjoint-VF-wave} 
	\langle \Box \widetilde{v} , {\mathcal{E}} p_\varrho \rangle_{Q_-}
	= \langle u_\varrho - \overline{u} , v \rangle_{L^2(Q)} \quad
	\mbox{for all} \; v\in \mathcal{H}_{0;0,}(Q). 
\end{equation}
The optimality system to be solved covers the forward (generalized) wave
equation \eqref{eq:generalized-VF-wave}, the adjoint backward (generalized)
wave equation \eqref{eq:adjoint-VF-wave}, and the gradient equation
\eqref{gradient equation p z}. When considering
$w_{z_\varrho} = A^{-1} z_\varrho = - \varrho^{-1} p_\varrho$
we can eliminate the control by $z_\varrho = -\varrho^{-1}Ap_\varrho$ to
end up with the system to find
$(u_\varrho,p_\varrho) \in \mathcal{H}_{0;0,}(Q) \times H_{0;,0}^{1,1}(Q)$
such that
\begin{equation}\label{eq:optimality-system_energy-reg}
	\varrho^{-1} \, \langle A p_\varrho , q \rangle_Q +
	\langle \Box \widetilde{u}_\varrho , {\mathcal{E}} q \rangle_{Q_-} = 0,
	\quad
	- \langle \Box \widetilde{v} , {\mathcal{E}} p_\varrho ) +
	\langle u_\varrho  , v \rangle_{L^2(Q)} =
	\langle \overline{u} , v \rangle_{L^2(Q)}
\end{equation}
is satisfied for all $(v,q) \in  
{\mathcal{H}}_{0;0,}(Q)\times H^{1,1}_{0;,0}(Q)$.

When the state $u_\varrho \in  \mathcal{H}_{0;0,}(Q)$ is known, and since
we are interested in the reconstruction of the control, we can compute 
$z_\varrho \in [H^{1,1}_{0;,0}(Q)]^*$ as unique solution of the
variational formulation
\begin{equation}\label{VF control}
	\langle z_\varrho , q \rangle_Q =
	\langle \Box \widetilde{u}_\varrho , {\mathcal{E}} q \rangle_{Q_-} 
	\quad \mbox{for all} \; q \in H^{1,1}_{0;,0}(Q).
\end{equation}

\section{Regularization error estimates}\label{Section Regularization}
We introduce $X:=\mathcal{H}_{0;0,}(Q)$ and $Y:=H_{0;,0}^{1,1}(Q)$ with
norms
\[
\| u \|_X = \| u \|_{{\mathcal{H}}_{0;0,}(Q)}, \quad
\| q \|_Y = \| q \|_{H^{1,1}_{0;,0}(Q)} = |q|_{H^1(Q)} ,
\]
and we can write the optimality system \eqref{eq:optimality-system_energy-reg}
as operator equation to find
$(u_\varrho,p_\varrho) \in X\times Y$ such that
\begin{equation}\label{eq:general-optimality-system}
	\begin{pmatrix}
		\varrho^{-1}A& B\\[1mm]
		-B^\ast & I
	\end{pmatrix}\begin{pmatrix}
		p_\varrho\\[1mm]
		u_\varrho
	\end{pmatrix} = \begin{pmatrix}
		0\\[1mm]
		\overline{u}
	\end{pmatrix}, 
\end{equation}
where $B : X \to Y^*$ is defined via
\[
\langle B v , q \rangle_Q = \langle \Box \widetilde{v} , {\mathcal{E}} q
\rangle_{Q_-} \quad \mbox{for all} \; (v,q) \in X\times Y.
\]
Note that, using \eqref{Boundedness Bilinear Form}, we have
\[
\| B u \|_{Y^*} = \sup\limits_{0 \neq q \in Y}
\frac{\langle Bu , q \rangle_Q}{\| q \|_Y} =
\sup\limits_{0 \neq q \in H^{1,1}_{0;,0}(Q)} \frac{\langle \Box \widetilde{u} ,
	{\mathcal{E}} q \rangle_{Q_-}}{\| q \|_{H^{1,1}_{0;,0}(Q)}} \leq
2 \, \| u \|_{{\mathcal{H}}_{0;0,}(Q)} = 2 \, \| u \|_X
\]
for all $u \in X$, i.e., $ B : X \to Y^*$ is bounded. 
Since $A$ is invertible, we can eliminate $p_\varrho = -\varrho A^{-1}Bu_\varrho$
to end up with the Schur complement equation to find
$u_\varrho \in X$ such that
\begin{equation}\label{eq:schur-complement} 
	\big[ I +\varrho B^* A^{-1}B\big]u_\varrho = \overline{u}
	\quad \mbox{in} \; X^*. 
\end{equation}

\begin{lemma}
	The operator $S := B^* A^{-1} B : X \to X^*$ is bounded and
	elliptic, i.e.,
	\[
	\| S u \|_{X^*} \leq 4 \, \| u \|_X, \quad
	\langle S u , u \rangle_Q \geq \| u \|_{X}^2
	\quad \mbox{for all} \; u \in X .
	\]
	Moreover, $ \| u \|_S := \langle Su,u \rangle_Q^{1/2}$, $u \in X$,
	defines an equivalent norm on $X$,
	\begin{equation}\label{norm equivalence X S}
		\| u \|_X \leq \| u \|_S \leq 2 \, \| u \|_X \quad \mbox{for all} \; u \in X.
	\end{equation}
\end{lemma}

\proof{
	The boundedness results from the boundedness of 
	$B:X \to Y^*$, and from the invertibility of $A : Y \to Y^*$, i.e.,
	for $ u \in X $ we have
	\[
	\| S u \|_{X^*} = \sup\limits_{0 \neq v \in X}
	\frac{\langle S u , v \rangle_Q}{\| v \|_X}
	= \sup\limits_{0 \neq v \in X}
	\frac{\langle A^{-1} B u , B v \rangle_Q}{\| v \|_X} \leq 4 \, \| u \|_X .
	\]
	Further, the inf-sup stability condition \eqref{inf-sup generalized}
	implies
	\[
	\| u \|_X =
	\| u \|_{\mathcal{H}_{0;0,}(Q)} \leq
	\sup_{0 \neq q \in H^{1,1}_{0;0,}(Q)}
	\frac{\langle \Box \widetilde{u},{\mathcal{E}} q\rangle_{Q_-}}
	{\|q\|_{H^{1,1}_{0;,0}(Q)}} =
	\sup_{0 \neq q \in Y}
	\frac{\langle Bu,q\rangle_Q}{\|q\|_Y} = \| Bu \|_{Y^*}
	\]
	for all $u\in X$. When introducing, for $u \in X$, the auxiliary variable
	$p_u = A^{-1} B u \in Y$, we first have
	\[
	\langle Su,u \rangle_Q =
	\langle A^{-1}Bu,Bu \rangle_Q =
	\langle p_u, Bu \rangle_Q =
	\langle p_u, A p_u\rangle_Q = \| p_u \|_Y^2 \, .
	\]
	Moreover, we have that 
	\[
	\| u \|_X \leq \| Bu \|_{Y^*} = \| Ap_u \|_{Y^*} \leq \| p_u \|_Y ,
	\]
	and we conclude
	\[
	\| u \|_X^2 \leq \| p_u \|_Y^2 =
	\langle Su,u \rangle_Q .
	\]
	This also shows that
	\[
	\|u\|_X^2 \leq \| u \|_S^2 = \langle Su,u\rangle _Q
	\leq \| S u\|_{X^*} \| u \|_X \leq 4 \, \|u\|^2_X \quad
	\mbox{for all} \; u \in X,
	\]
	which gives the desired equivalence of norms.
}

\noindent
The variational formulation of the Schur complement equation
\eqref{eq:schur-complement} is to find $ u_\varrho \in X $ such that
\begin{equation}\label{eq:schur-complement VF}
	\varrho \, \langle S u_\varrho , v \rangle_Q +
	\langle u_\varrho , v \rangle_{L^2(Q)} =
	\langle \overline{u} , v \rangle_{L^2(Q)} \quad
	\mbox{for all} \; v \in X .
\end{equation}
Unique solvability of \eqref{eq:schur-complement VF} immediately
follows from the properties of $S$ for all $\varrho \in {\mathbb{R}}_+$.
In particular for $v = u_\varrho$ this gives
\[
\varrho \, \| u_\varrho \|_S^2 + \| u_\varrho \|^2_{L^2(Q)} =
\langle \overline{u} , u_\varrho \rangle_{L^2(Q)} \leq
\| \overline{u} \|_{L^2(Q)} \| u_\varrho \|_{L^2(Q)},
\]
and hence,
\begin{equation}\label{eq:H1L2estimate-continuous}    
	\| u_\varrho \|_{L^2(Q)} \leq \| \overline{u} \|_{L^2(Q)}, \quad
	\sqrt{\varrho} \, \| u_\varrho \|_S \leq \| \overline{u} \|_{L^2(Q)} 
\end{equation}
follow. As in \cite[Lemma 2.3]{UL:LangerSteinbachYang:2022a}
we can prove the following regularization error  estimates.

\begin{theorem}\label{thm:error-estimates}
	For $\overline{u} \in L^2(Q)$ let $u_\varrho \in X$ be the unique solution
	of the variational formulation \eqref{eq:schur-complement VF} where
	$\varrho \in {\mathbb{R}}_+$. Then the
	following estimate holds true
	\begin{equation}\label{eq:l2l2estimate-continuous} 
		\| u_\varrho-\overline{u} \|_{L^2(Q)} \leq
		\| \overline{u} \|_{L^2(Q)}.
	\end{equation}
	Moreover, for $ \overline{u} \in X$ we have
	\begin{equation}\label{eq:l2H1estimate-continuous}
		\| u_\varrho - \overline{u} \|_{L^2(Q)} \leq \sqrt{\varrho} \,
		\| \overline{u} \|_S \, ,
	\end{equation}
	as well as
	\begin{equation}\label{eq:H1H1estimate-continuous}
		\| u_\varrho - \overline{u} \|_S \leq \| \overline{u} \|_S \, .
	\end{equation}
	If in addition $\overline{u} \in X$ is such that
	$S \overline{u} \in L^2(Q)$ is satisfied, then
	\begin{equation} \label{eq:l2Sestimate-continuous}
		\| u_\varrho - \overline{u} \|_{L^2(Q)} \leq
		\varrho \, \| S \overline{u} \|_{L^2(Q)},
	\end{equation}
	and 
	\begin{equation}\label{eq:SSestimate-continuous}
		\| u_\varrho - \overline{u} \|_S \leq \sqrt{\varrho} \,
		\| S \overline{u} \|_{L^2(Q)}. 
	\end{equation}
\end{theorem}

\proof{
	Let us first consider the case $\overline u \in L^2(Q)$.
	Then, when choosing $v=u_\varrho \in X$ within the variational formulation 
	\eqref{eq:schur-complement VF}, this gives
	\[
	\varrho \, \langle S u_\varrho , u_\varrho \rangle_Q =
	\langle \overline{u} - u_\varrho , u_\varrho \rangle _{L^2(Q)} =
	- \langle \overline{u} - u_\varrho , \overline{u}-u_\varrho
	\rangle_{L^2(Q)} +
	\langle \overline{u} - u_\varrho,\overline{u} \rangle _{L^2(Q)} ,
	\]
	i.e.,
	\[
	\| u_\varrho - \overline{u} \|_{L^2(Q)}^2 + \varrho \, \| u_\varrho \|^2_S =
	\langle \overline{u} - u_\varrho , u_\varrho \rangle_{L^2(Q)} \leq
	\| u_\varrho - \overline{u} \|_{L^2(Q)} \| \overline{u} \|_{L^2(Q)} ,
	\]
	and \eqref{eq:l2l2estimate-continuous} follows.
	
	For $\overline{u} \in X$ we can consider the variational formulation
	\eqref{eq:schur-complement VF} for $v=\overline{u} - u_\varrho \in X$ to
	obtain
	\begin{eqnarray*}
		\| \overline{u} - u_\varrho \|^2_{L^2(Q)}
		& = & \langle \overline{u} - u_\varrho , \overline{u} - u_\varrho
		\rangle_{L^2(Q)} \, = \,
		\varrho \, \langle S u_\varrho , \overline{u} - u_\varrho
		\rangle_Q \\
		& = & \varrho \, \langle S \overline{u} , \overline{u} - u_\varrho
		\rangle_Q - \varrho \, \langle S (\overline{u} - u_\varrho),
		\overline{u} - u_\varrho \rangle_Q ,
	\end{eqnarray*}
	i.e.,
	\[
	\| u_\varrho - \overline{u} \|^2_{L^2(Q)} +
	\varrho \, \| u_\varrho - \overline{u} \|^2_S \leq
	\varrho \, \langle S \overline{u} , \overline{u} - u_\varrho \rangle_Q
	\leq \varrho \, \| \overline{u} \|_S \| u_\varrho - \overline{u} \|_S ,
	\]
	and hence, \eqref{eq:H1H1estimate-continuous} and
	\eqref{eq:l2H1estimate-continuous} follow.
	
	If $\overline{u} \in X$ is such that $S \overline{u} \in L^2(Q)$ is
	satisfied, we also have
	\[
	\| u_\varrho - \overline{u} \|^2_{L^2(Q)} +
	\varrho \, \| u_\varrho - \overline{u} \|^2_S \leq
	\varrho \, \langle S \overline{u} , \overline{u} - u_\varrho \rangle_Q
	\leq \varrho \, \| S \overline{u} \|_{L^2(Q)}
	\| u_\varrho - \overline{u} \|_{L^2(Q)} ,
	\]
	from which \eqref{eq:l2Sestimate-continuous} and
	\eqref{eq:SSestimate-continuous} follow.
}

\begin{corollary}\label{corollary regularization}
	For $\overline{u} \in H^{1,1}_{0;0,}(Q) \subset X = {\mathcal{H}}_{0;0,}(Q)$
	we conclude from \eqref{eq:l2H1estimate-continuous},
	\eqref{norm equivalence X S}, and \eqref{norm X H1},
	\[
	\| u_\varrho - \overline{u} \|_{L^2(Q)} \leq \sqrt{\varrho} \,
	\| \overline{u} \|_S \, \leq 2 \, \sqrt{\varrho} \,
	\| \overline{u} \|_X \leq 2 \, \sqrt{\varrho} \,
	\| \overline{u} \|_{H^{1,1}_{0;0,}(Q)} ,
	\]
	and using a space interpolation argument,
	see, e.g., \cite{Adams:1975, LionsMagenesVol1:1968, McLean:2000},
	this gives
	\begin{equation}\label{regularization error Hs}
		\| u_\varrho - \overline{u} \|_{L^2(Q)} \leq c \, \varrho^{s/2} \,
		\| \overline{u} \|_{H^{s,s}_{0;0,}(Q)}
	\end{equation}
	when assuming $\overline{u} \in H^{s,s}_{0;0,}(Q) :=
	[L^2(Q),H^{1,1}_{0;0,}(Q)]_s$ for some $s \in [0,1]$, and where the
	positive constant $c$ is independent of $\varrho$.
	
	Next we consider $\overline{u} \in H^{1,1}_{0;0,}(Q) \cap H^2(Q)$ and
	assume that $\overline{u}$ is such that
	$A^{-1} B \overline{u} \in H^2(Q)$.  Note that $A$ is
	related to the space-time Laplacian, but with mixed Dirichlet and
	Neumann boundary conditions which may reduce the regularity of its
	solution. The application
	of the adjoint wave operator then finally gives
	$S \overline{u} = B^* A^{-1} B \overline{u} \in L^2(Q)$.
	Then the error estimate
	\eqref{eq:l2Sestimate-continuous} implies
	\[
	\| u_\varrho - \overline{u} \|_{L^2(Q)} \leq c \, \varrho \,
	\| \overline{u} \|_{H^2(Q)},
	\]
	and using an interpolation argument finally gives
	\begin{equation}
		\| u_\varrho - \overline{u} \|_{L^2(Q)} \leq c \, \varrho^{s/2} \,
		\| \overline{u} \|_{H^s(Q)}
	\end{equation}
	when assuming $ \overline{u} \in H^{1,1}_{0;0,}(Q) \cap H^s(Q)$ for
	some $s \in (1,2]$.
\end{corollary}

\section{Space-time finite element methods}
For the Galerkin discretization of the Schur complement variational
formulation \eqref{eq:schur-complement VF} we introduce the conforming
finite element space
$X_h := S_h^1(Q)\cap \mathcal{H}_{0;0,}(Q) =
\mbox{span} \{ \varphi_k \}_{k=1}^{M_X} \subset X$ of
piecewise linear and continuous basis functions $\varphi_k$ which are
defined with respect to some admissible globally quasi-uniform decomposition $\mathcal{T}_h=\{\tau_\ell\}_{\ell=1}^N$
of the space-time domain $Q$ into shape-regular simplicial finite elements $\tau_\ell$
of mesh size $h_\ell$, see, e.g., \cite{BrennerScott:2008}.
Then the finite element approximation of \eqref{eq:schur-complement VF}
is to find $u_{\varrho h} \in X_h$ such that
\begin{equation}\label{Schur VF FEM}
	\varrho \, \langle S u_{\varrho h} , v_h \rangle_Q +
	\langle u_{\varrho h} , v_h \rangle_{L^2(Q)} =
	\langle \overline{u} , v_h \rangle_{L^2(Q)}
\end{equation}
ist satisfied for all $v_h \in X_h$. Using standard arguments, we conclude
unique solvability of \eqref{Schur VF FEM}, and the following
Cea type a priori error estimate
\begin{equation}\label{FEM Cea}
	\varrho \, \| u_\varrho - u_{\varrho h} \|_S^2 +
	\| u_\varrho - u_{\varrho h} \|_{L^2(Q)}^2
	\leq
	\inf\limits_{v_h \in X_h} \Big[ \varrho \, \| u_\varrho - v_h \|_S^2 +
	\| u_\varrho - v_h \|_{L^2(Q)}^2 \Big] \, .
\end{equation}

\begin{theorem}
	Assume $\overline{u} \in [L^2(Q), H^{1,1}_{0;0,}(Q)]_s$ for $s \in [0,1]$
	or $\overline{u} \in H^{1,1}_{0;0,}(Q) \cap H^s(Q)$ for $s \in (1,2]$.
	For the unique solution $u_{\varrho h} \in X_h$ of
	\eqref{Schur VF FEM} there holds the finite element error estimate
	\begin{equation}\label{FEM error L2}
		\| u_{\varrho h} - \overline{u} \|_{L^2(Q)} \leq c \, h^s \,
		\|\overline{u} \|_{H^s(Q)},
	\end{equation}
	provided that $\varrho = h^2$. For
	$\overline{u} \in H^{1,1}_{0;0,}(Q) \cap H^s(Q)$ and $s \in [1,2]$ we
	also have the error estimate
	\begin{equation}\label{FEM error S}
		\| u_{\varrho h} - \overline{u} \|_S \leq c \, h^{s-1} \,
		\|\overline{u} \|_{H^s(Q)}.
	\end{equation}
\end{theorem}

\proof{
	We first consider the error estimate \eqref{FEM Cea} for the
	particular function $v_h \equiv 0$, and using
	\eqref{eq:H1L2estimate-continuous} this gives
	\[
	\| u_\varrho - u_{\varrho h} \|^2_{L^2(Q)} \leq
	\varrho \, \| u_\varrho \|_S^2 + \| u_\varrho \|_{L^2(Q)}^2 \leq
	2 \, \| \overline{u} \|^2_{L^2(Q)} \, .
	\]
	Hence we conclude
	\begin{equation}\label{FEM error L2 L2}
		\| u_{\varrho h} - \overline{u} \|_{L^2(Q)} \leq
		\| u_\varrho - \overline{u} \|_{L^2(Q)} +
		\| u_\varrho - u_{\varrho h} \|_{L^2(Q)} \leq (1+\sqrt{2}) \,
		\| \overline{u} \|_{L^2(Q)} \, .
	\end{equation}
	We now assume $\overline{u} \in H^{1,1}_{0;0,}(Q) \subset X$,
	and from \eqref{FEM Cea} we obtain, using the triangle inequality,
	\eqref{eq:H1H1estimate-continuous} and
	\eqref{eq:l2H1estimate-continuous}, the inclusion
	$H^{1,1}_{0;0,}(Q) \subset X$, and standard approximation properties
	of piecewise linear finite element functions, e.g.,
	Scott--Zhang interpolation \cite{BrennerScott:2008},
	\begin{eqnarray*}
		\varrho \, \| u_\varrho - u_{\varrho h} \|_S^2 +
		\| u_\varrho - u_{\varrho h} \|^2_{L^2(Q)}
		& \leq & \inf\limits_{v_h \in X_h} \Big[
		\varrho \, \| u_\varrho - v_h \|^2_S +
		\| u_\varrho - v_h \|^2_{L^2(Q)} \Big] \\
		&& \hspace*{-5cm} \leq 2 \left[
		\varrho \, \| u_\varrho - \overline{u} \|^2_S +
		\| u_\varrho - \overline{u} \|_{L^2(Q)}^2 +
		\inf\limits_{v_h \in X_h} \Big[ 
		\varrho \, \| \overline{u} - v_h \|_S^2 +
		\| \overline{u} - v_h \|^2_{L^2(Q)}
		\Big]
		\right] \\
		&& \hspace*{-5cm} \leq 2 \left[ 2 \,
		\varrho \, \| \overline{u} \|^2_S +
		\inf\limits_{v_h \in X_h} \Big[ c \,
		\varrho \, \| \overline{u} - v_h \|_{H^{1,1}_{0;0,}(Q)}^2 +
		\| \overline{u} - v_h \|^2_{L^2(Q)}
		\Big]
		\right] \\
		&& \hspace*{-5cm} \leq c \, \Big[ \varrho + h^2 \Big] \,
		\| \overline{u} \|_{H^1(Q)}^2 \, .
	\end{eqnarray*}
	In particular for $\varrho = h^2$ this gives
	\[
	h^2 \, \| u_\varrho - u_{\varrho h} \|_S^2 +
	\| u_\varrho - u_{\varrho h} \|^2_{L^2(Q)} \leq c \, h^2 \,
	\| \overline{u} \|_{H^1(Q)}^2 \, .
	\]
	Hence, using the triangle inequality and
	\eqref{regularization error Hs},
	\begin{equation}\label{FEM error L2 H1}
		\| u_{\varrho h} - \overline{u} \|_{L^2(Q)} \leq
		\| u_{\varrho h} - u_\varrho \|_{L^2(Q)} +
		\| u_\varrho - \overline{u} \|_{L^2(Q)} \leq
		c \, h \,
		\| \overline{u} \|_{H^1(Q)} 
	\end{equation}
	follows, while with \eqref{eq:H1H1estimate-continuous} we obtain
	\begin{equation}\label{FEM error H1 H1}
		\| u_{\varrho h} - \overline{u} \|_S \leq
		\| u_{\varrho h} - u_\varrho \|_S +
		\| u_\varrho - \overline{u} \|_S \leq
		c \, \| \overline{u} \|_{H^1(Q)} .
	\end{equation}
	For $ \overline{u} \in H^{1,1}_{0;0,}(Q) \cap H^2(Q) \subset X$,
	using \eqref{eq:SSestimate-continuous} and
	\eqref{eq:l2Sestimate-continuous}, we can prove in the same way
	\[
	\varrho \, \| u_\varrho - u_{\varrho h} \|_S^2 +
	\| u_\varrho - u_{\varrho h} \|^2_{L^2(Q)} \leq c \, \Big[
	\varrho^2 + \varrho \, h^2 + h^4 \Big] \, \| \overline{u} \|_{H^2(Q)}^2 =
	c \, h^4 \, \| \overline{u} \|^2_{H^2(Q)},
	\]
	provided that $\varrho=h^2$. Now, using
	\eqref{eq:l2Sestimate-continuous},
	\eqref{eq:SSestimate-continuous} and
	Corollary \ref{corollary regularization},
	we obtain
	\[
	\| u_{\varrho h} - \overline{u} \|_{L^2(Q)} \leq c \, h^2 \,
	\| \overline{u} \|_{H^2(Q)},
	\]
	and
	\[
	\| u_{\varrho h} - \overline{u} \|_S \leq c \, h \,
	\| \overline{u} \|_{H^2(Q)} .
	\]
	The general estimates for $s \in (0,1]$ and $s \in (1,2)$ now follow
	again from a space interpolation argument.
}

\begin{corollary}
	As already given in the previous proof, there hold the error estimates
	\[
	\varrho \, \| u_\varrho - u_{\varrho h} \|_S^2 +
	\| u_\varrho - u_{\varrho h} \|^2_{L^2(Q)}
	\leq c \, \Big[ \varrho + h^2 \Big] \,
	\| \overline{u} \|_{H^1(Q)}^2 
	\]
	when assuming $\overline{u} \in H^{1,1}_{0;0,}(Q)$, and
	\[
	\varrho \, \| u_\varrho - u_{\varrho h} \|_S^2 +
	\| u_\varrho - u_{\varrho h} \|^2_{L^2(Q)} \leq c \, \Big[
	\varrho^2 + \varrho \, h^2 + h^4 \Big] \, \| \overline{u} \|_{H^2(Q)}^2
	\]
	when assuming $\overline{u} \in H^{1,1}_{0;0,}(Q) \cap H^2(Q)$.
\end{corollary}

\noindent
Next we are going to define a computable approximation of
$Su = B^* A^{-1} Bu$. For $u \in X$, let
$p_u = A^{-1} B u \in Y$ be the unique solution of the variational
formulation
\[
\langle Ap_u, q \rangle_Q = \langle Bu,q\rangle _Q
\quad \mbox{for all} \; q \in Y ,
\]
and hence, $S u = B^* p_u$. Let
$Y_h:= S_h^1(Q)\cap H_{0;,0}^{1,1}(Q) = \mbox{span} \{ \psi_i \}_{i=1}^{M_Y}$
be a second finite element space of piecewise linear continuous basis functions,
which, for simplicity, are defined with respect to the same decomposition
of the space-time domain $Q$ into finite elements as $X_h$.
Let now $p_{uh}\in Y_h$ solve
\[
\langle Ap_{uh} ,q_h \rangle_Q =
\langle Bu,q_h\rangle_Q \quad \mbox{for all} \; q_h \in Y_h ,
\]
and define $\widetilde{S}u := B^* p_{uh}$, where
$\widetilde{S} : X \to X^*$ is bounded due to the properties of
$A:Y \to Y^*$ and $B:X \to Y^*$, respectively. Instead of
\eqref{Schur VF FEM}, we now consider the perturbed variational formulation
to find $\widetilde{u}_{\varrho h} \in X_h$ such that
\begin{equation}\label{Schur VF FEM pert}
	\varrho \, \langle \widetilde{S} \widetilde{u}_{\varrho h} , v_h \rangle_Q
	+ \langle \widetilde{u}_{\varrho h} , v_h \rangle_{L^2(Q)} =
	\langle \overline{u} , v_h \rangle_{L^2(Q)}
\end{equation}
is satisfied for all $v_h \in X_h$. Unique solvability of
\eqref{Schur VF FEM pert} follows since the matrix realization of
$\widetilde{S}$ is positive semi-definite, while the mass matrix,
which is related to the inner product in $L^2(Q)$, is positive definite.

\begin{lemma}\label{lem:strang-lemma}
	Let $u_{\varrho h}, \widetilde{u}_{\varrho h} \in X_h$ be the
	unique solutions of the variational formulations \eqref{Schur VF FEM}
	and \eqref{Schur VF FEM pert}, respectively. For $\overline{u} \in L^2(Q)$
	there holds the error estimate
	\begin{equation}\label{FEM pert L2}
		\| \widetilde{u}_{\varrho h} - \overline{u} \|_{L^2(Q)} \leq
		\| \overline{u} \|_{L^2(Q)} .
	\end{equation}
	For $\overline{u} \in H^{1,1}_{0;0,}(Q) \cap H^2(Q) \subset X$ we have 
	\begin{equation}\label{FEM pert H2}
		\| \widetilde{u}_{\varrho h} - \overline{u} \|_{L^2(Q)} \leq c \, h^2 \,
		\| \overline{u} \|_{H^2(Q)} .
	\end{equation}
\end{lemma}
\proof{
	The estimate \eqref{FEM pert L2} follows when considering the
	perturbed variational formulation \eqref{Schur VF FEM pert} for
	$v_h = \widetilde{u}_{\varrho h}$, i.e.,
	\[
	\varrho \, \langle \widetilde{S} \widetilde{u}_{\varrho h} ,
	\widetilde{u}_{\varrho h} \rangle_Q
	=
	\langle \overline{u} - \widetilde{u}_{\varrho h},
	\widetilde{u}_{\varrho h} \rangle_{L^2(Q)}
	=
	- \langle \overline{u} - \widetilde{u}_{\varrho h},
	\overline{u} - \widetilde{u}_{\varrho h} \rangle_{L^2(Q)}
	+ \langle \overline{u} - \widetilde{u}_{\varrho h},
	\overline{u} \rangle_{L^2(Q)},
	\]
	and
	\[
	\varrho \, \langle \widetilde{S} \widetilde{u}_{\varrho h} ,
	\widetilde{u}_{\varrho h} \rangle_Q
	+ \| \overline{u} - \widetilde{u}_{\varrho h} \|^2_{L^2(Q)}
	= \langle \overline{u} - \widetilde{u}_{\varrho h},
	\overline{u} \rangle_{L^2(Q)} \leq \| \overline{u} -
	\widetilde{u}_{\varrho h} \|_{L^2(Q)} \| \overline{u} \|_{L^2(Q)} .
	\]
	When subtracting the perturbed variational formulation
	\eqref{Schur VF FEM pert} from \eqref{Schur VF FEM}, this gives
	\[
	\varrho \, \langle S u_{\varrho h} -
	\widetilde{S} \widetilde{u}_{\varrho h} , v_h \rangle_Q +
	\langle u_{\varrho h} - \widetilde{u}_{\varrho h} , v_h
	\rangle_{L^2(Q)} = 0 \quad \mbox{for all} \; v_h \in X_h ,
	\]
	i.e.,
	\[
	\varrho \, \langle (S-\widetilde{S}) u_{\varrho h}, v_h \rangle_Q +
	\langle u_{\varrho h} - \widetilde{u}_{\varrho h} , v_h
	\rangle_{L^2(Q)} =
	\varrho \, \langle \widetilde{S} (\widetilde{u}_{\varrho h} - u_{\varrho h}),
	v_h \rangle_Q
	\quad \mbox{for all} \; v_h \in X_h .
	\]
	In particular for $v_h = \widetilde{u}_{\varrho h} - u_{\varrho h}$ we
	further conclude
	\begin{eqnarray*}
		0 & \leq & \varrho \,
		\langle \widetilde{S} (\widetilde{u}_{\varrho h} -
		u_{\varrho h}),
		\widetilde{u}_{\varrho h} - u_{\varrho h} \rangle_Q \\
		& = & \varrho \, \langle (S-\widetilde{S}) u_{\varrho h},
		\widetilde{u}_{\varrho h} - u_{\varrho h} \rangle_Q +
		\langle u_{\varrho h} - \widetilde{u}_{\varrho h} ,
		\widetilde{u}_{\varrho h} - u_{\varrho h} \rangle_{L^2(Q)},
	\end{eqnarray*}
	i.e., using an inverse inequality in $X_h$,
	\begin{eqnarray*}
		\| \widetilde{u}_{\varrho h} - u_{\varrho h} \|^2_{L^2(Q)}
		& \leq & \varrho \, \langle (S-\widetilde{S}) u_{\varrho h},
		\widetilde{u}_{\varrho h} - u_{\varrho h} \rangle_Q \\
		& = & \varrho \, \langle B^* (p_{u_{\varrho h}} - p_{u_{\varrho h}h}),
		\widetilde{u}_{\varrho h} - u_{\varrho h} \rangle_Q \\
		& = & \varrho \, \langle p_{u_{\varrho h}} - p_{u_{\varrho h}h},
		B(\widetilde{u}_{\varrho h} - u_{\varrho h}) \rangle_Q \\
		& = & \varrho \, \langle \Box (\widetilde{u}_{\varrho h} -
		u_{\varrho h}) ,
		{\mathcal{E}}(p_{u_{\varrho h}} - p_{u_{\varrho h}h})
		\rangle_{Q_-} \\
		& = & \varrho \, b( \widetilde{u}_{\varrho h} - u_{\varrho h},
		p_{u_{\varrho h}} - p_{u_{\varrho h}h}) \\
		& \leq & \varrho \, | \widetilde{u}_{\varrho h} - u_{\varrho h} |_{H^1(Q)}
		| p_{u_{\varrho h}} - p_{u_{\varrho h}h} |_{H^1(Q)} \\
		& \leq & c \, \varrho \, h^{-1} \,
		\| \widetilde{u}_{\varrho h} - u_{\varrho h} \|_{L^2(Q)}
		| p_{u_{\varrho h}} - p_{u_{\varrho h}h} |_{H^1(Q)}.
	\end{eqnarray*}
	Hence, using $\varrho = h^2$ and the triangle inequality, this gives
	\begin{eqnarray*}
		\| \widetilde{u}_{\varrho h} - u_{\varrho h} \|_{L^2(Q)}
		& \leq & c \, h \, | p_{u_{\varrho h}} - p_{u_{\varrho h}h} |_{H^1(Q)} \\
		& \leq & c \, h \, \Big[
		| p_{u_{\varrho h}} - p_{\overline{u}} |_{H^1(Q)}
		+ | p_{\overline{u}} - p_{\overline{u} h} |_{H^1(Q)}
		+ | p_{\overline{u}h} - p_{u_{\varrho h}h} |_{H^1(Q)}
		\Big] \, .
	\end{eqnarray*}
	For the first term we further have
	\begin{eqnarray*}
		| p_{u_{\varrho h}} - p_{\overline{u}} |_{H^1(Q)}^2
		& = & \langle A (  p_{u_{\varrho h}} - p_{\overline{u}} ),
		p_{u_{\varrho h}} - p_{\overline{u}} \rangle_Q \\
		& = & \langle B (  u_{\varrho h} - \overline{u} ),
		p_{u_{\varrho h}} - p_{\overline{u}} \rangle_Q \\
		& \leq & \| B (  u_{\varrho h} - \overline{u} ) \|_{Y^*}
		\| p_{u_{\varrho h}} - p_{\overline{u}} \|_Y \\
		& \leq & 2 \, \| u_{\varrho h} - \overline{u} \|_X
		\| p_{u_{\varrho h}} - p_{\overline{u}} \|_Y,
	\end{eqnarray*}
	i.e.,
	\[
	| p_{u_{\varrho h}} - p_{\overline{u}} |_{H^1(Q)} \leq
	2 \, \| u_{\varrho h} - \overline{u} \|_X \leq
	2 \, \| u_{\varrho h} - \overline{u} \|_S \leq c \, h \,
	\| \overline{u} \|_{H^2(Q)} .
	\]
	Following the same lines we can also estimate the third term by 
	\[
	| p_{\overline{u}h} - p_{u_{\varrho h}h} |_{H^1(Q)}
	\leq c \, h \, \| \overline{u} \|_{H^2(Q)} .
	\]
	To estimate the second term, let us first recall that
	$p_{\overline u}\in Y = H_{0;,0}^{1,1}(Q)$ solves 
	\[
	\langle Ap_{\overline{u}}, q \rangle_Q =
	\langle B \overline{u} , q \rangle_Q \quad
	\mbox{for all} \; q \in Y,
	\]
	while $p_{\overline{u}h}\in Y_h$ solves 
	\[
	\langle A p_{\overline{u}h}, q_h \rangle_Q =
	\langle B \overline{u}, q_h \rangle_Q \quad
	\mbox{for all} \; q_h \in Y_h.
	\]
	Thus, we conclude the Galerkin orthogonality
	\[
	\langle A( p_{\overline{u}} - p_{\overline{u}h}) ,q_h \rangle _Q =0
	\quad \mbox{for all} \; q_h \in Y_h,
	\]
	and Cea's lemma,
	\[
	| p_{\overline{u}} - p_{\overline{u}h}|_{H^1(Q)} \leq
	\inf\limits_{q_h \in Y_h} | p_{\overline{u}} - q_h|_{H^1(Q)} \leq
	c \, h \, |p_{\overline{u}}|_{H^2(Q)} ,
	\]
	when assuming $p_{\overline{u}} = A^{-1} B \overline{u} \in H^2(Q)$.
	Indeed, for a convex space-time domain $Q$ we have
	\[
	|p_{\overline{u}}|_{H^2(Q)} \leq c \,
	\| A p_{\overline{u}} \|_{L^2(Q)} = c \,
	\| B \overline{u} \|_{L^2(Q)} \leq c \, \| \overline{u} \|_{H^2(Q)} .
	\]
	This concludes the proof.
}

\begin{remark}
	In the proof of Lemma~\ref{lem:strang-lemma} we have used an inverse
	inequality which in general assumes a globally quasi-uniform
	finite element mesh. However, when using a variable regularization
	function $\varrho(x,t) = h_\ell^2$ for $(x,t) \in \tau_\ell$, it is
	sufficient to use the inverse inequality locally, allowing
	adaptively refined and locally quasi-uniform finite element meshes.
	For a related approach for a distributed optimal control problem with
	variable regularization
	subject to the Poisson equation, see \cite{LLSY:2022variable}.
\end{remark}

\begin{corollary}
	When using a space interpolation argument, from \eqref{FEM pert L2}
	and \eqref{FEM pert H2} we now conclude the final error estimate
	\begin{equation}\label{FEM pert Hs}
		\| \widetilde{u}_{\varrho h} - \overline{u} \|_{L^2(Q)}
		\leq c \, h^s \ \| \overline{u} \|_{H^s(Q)} 
	\end{equation}
	when assuming $\overline{u} \in [L^2(Q),H^{1,1}_{0;0,}(Q)]_s$ for
	$s \in [0,1]$ or $\overline{u} \in H^{1,1}_{0;0,}(Q) \cap H^s(Q)$
	for $s \in (1,2]$.
\end{corollary}

\noindent
When the approximate state $\widetilde{u}_{\varrho h}$ is known,
as in \eqref{VF control} we can compute the associate control
$\widetilde{z}_\varrho = B \widetilde{u}_{\varrho h}
\in [H^{1,1}_{0;,0}(Q)]^*$ as unique solution
of the variational formulation
\[
\langle \widetilde{z}_\varrho , q \rangle_Q =
\langle \Box \widetilde{u}_{\varrho h} , {\mathcal{E}}q \rangle_{Q_-} =
\langle B \widetilde{u}_{\varrho h} , q \rangle_Q
\quad \mbox{for all} \; q \in H^{1,1}_{0;,0}(Q) .
\]
With this we conclude that
$\widetilde{z}_\varrho$ is the minimizer of the functional
\[
{\mathcal{F}}(z) :=
\frac{1}{2} \, \| z - B \widetilde{u}_{\varrho h}
\|^2_{[H^{1,1}_{0;,0}(Q)]^*} = \frac{1}{2} \,
\langle A^{-1} (z - B \widetilde{u}_{\varrho h}),
z - B \widetilde{u}_{\varrho h} \rangle_Q,
\]
i.e., $\widetilde{z}_{\varrho} \in [H^{1,1}_{0;,0}(Q)]^*$ is the unique
solution of the gradient equation
\[
A^{-1} (\widetilde{z}_{\varrho} - B \widetilde{u}_{\varrho h}) = 0.
\]
This is equivalent to the coupled system to find
$(\psi,\widetilde{z}_\varrho) \in H^{1,1}_{0;,0}(Q) \times
[H^{1,1}_{0;,0}(Q)]^*$ such that
\begin{equation}\label{eq:VF_control-continuous}
	A \psi + \widetilde{z}_{\varrho} = B \widetilde{u}_{\varrho h},
	\quad \psi = 0.
\end{equation}
Let $Z_H \subset [H^{1,1}_{0;,0}(Q)]^*$ be a suitable finite element
space, then we consider the Galerkin variational formulation to find
$(\psi_h,\widetilde{z}_{\varrho H}) \in Y_h \times Z_H$ such that
\begin{equation}\label{eq:VF_control-discrete}
	\langle A \psi_h , \phi_h \rangle_{L^2(Q)} +
	\langle \widetilde{z}_{\varrho H} , \phi_h \rangle_{L^2(Q)} =
	\langle B \widetilde{u}_{\varrho h} , \phi_h \rangle_{L^2(Q)}, \quad
	\langle \psi_h , \eta_H \rangle_{L^2(Q)} = 0 
\end{equation}
is satisfied for all $(\phi_h,\eta_H) \in Y_h \times Z_H$.
Unique solvability of \eqref{eq:VF_control-discrete} follows when
the discrete inf-sup stability condition
\[
c_S \, \| z_H \|_{[H^{1,1}_{0;,0}(Q)]^*} \leq
\sup\limits_{0 \neq\phi_h \in Y_h}
\frac{\langle z_H ,\phi_h \rangle_{L^2(Q)}}{\| \phi_h \|_{H^{1,1}_{0;,0}(Q)}}
\quad \mbox{for all} \; z_H \in Z_H
\]
is satisfied, i.e., when $Y_h$ is defined with respect to a space-time
finite element mesh size $h$ which is sufficiently small compared to the
mesh size $H$ of
$Z_H$. From a practical point of view it is sufficient to consider
one additional refinement when defining first $Z_H$, and afterwards $Y_h$,
i.e., $h=H/2$.
As in mixed finite element methods and using the Strang lemma we
can then derive related error estimates for the Galerkin solution
$\widetilde{z}_{\varrho H}$.

\section{Numerical results}
The perturbed variational formulation \eqref{Schur VF FEM pert}
corresponds to the Ga\-ler\-kin discretization of the coupled
variational formulation \eqref{eq:optimality-system_energy-reg}.
With the finite element spaces
\[
X_h:=S_h^1(\mathcal{T}_h)\cap H_{0;0,}^{1,1}(Q)=
\text{span}\{\varphi_k\}_{k=1}^{M_X}
\]
and
\[
Y_h:=S_h^1(\mathcal{T}_h)\cap
H_{0;,0}^{1,1}(Q)=\text{span}\{\psi_i\}_{i=1}^{M_Y}
\]
as already used in Section 4, the equivalent
linear system of algebraic equations reads
\begin{equation}\label{eq:matrix-form-energ-reg}
	\begin{pmatrix}
		\varrho^{-1} A_h & B_h\\
		-B_h^\top & M_h
	\end{pmatrix}
	\begin{pmatrix}
		\underline p\\\underline u
	\end{pmatrix}=\begin{pmatrix}
		\underline 0\\\underline{f}
	\end{pmatrix},
\end{equation}
where the system matrix is positive definite but skew-symmetric, and 
where the matrix entries are given as, for $k,\ell=1,\ldots,M_X$,
$i,j=1,\ldots,M_Y$,
\begin{eqnarray*}
	A_h[j,i]
	& = & \langle \nabla_{(x,t)}\psi_i,\nabla_{(x,t)}\psi_j\rangle_{L^2(Q)}, \\
	M_h[\ell,k]
	& = & \langle\varphi_k,\varphi_\ell\rangle_{L^2(Q)}, \\
	B_h[j,k]
	& = & - \langle \partial_t \varphi_k , \partial_t \psi_j \rangle_{L^2(Q)} +
	\langle \nabla_x \varphi_k , \nabla_x \psi_j \rangle_{L^2(Q)},
\end{eqnarray*}
and with the load vector
\begin{align*}
	f_\ell = \langle \overline u,\varphi_\ell \rangle_{L^2(Q)}.     
\end{align*}
In addition to the energy regularization
we will also consider the control
$z_\varrho\in L^2(Q)$, where the arising matrix system is given as 
\begin{equation}\label{eq:matrix-form-l2reg}
	\begin{pmatrix}
		\varrho^{-1} \overline M_h & B_h\\
		-B_h^\top & M_h
	\end{pmatrix}
	\begin{pmatrix}
		\underline p\\\underline u
	\end{pmatrix}=\begin{pmatrix}
		\underline 0\\\underline f
	\end{pmatrix},
\end{equation}
with the related mass matrix
\[
\overline M_h[j,i]:=\langle \psi_i,\psi_j\rangle_{L^2(Q)} \quad
\mbox{for} \; i,j=1,\ldots,M_Y.
\]
A similar analysis as for the energy regularization shows, that in this case,
the optimal choice for the relaxation parameter is $\varrho = h^4$, see
also \cite{LangerLoescherSteinbachYang:2022a} in the case of a distributed
optimal control problem for the Poisson equation.

\subsection{Uniform refinement}
\begin{figure}[htbp!]
	\centering
	\begin{subfigure}[b]{0.45\textwidth}
		\centering
		\begin{tikzpicture}[scale=0.7]
			\begin{axis}[
				xmode = log,
				ymode = log,
				xlabel={Number $N$
					of space-time finite elements},
				ylabel=$\|\widetilde u_{i,\varrho h}-\overline u_{i}\|_{L^2(Q)}$,
				legend style={font=\tiny}, legend pos = south west]
				\addplot table [col sep=&, y=err1, x=N]{tab-examples2D.dat};
				\addlegendentry{$\|\widetilde u_{1,\varrho h}-\overline u_{1}\|_{L^2(Q)}$}
				\addplot table [col sep=&, y=err2, x=N]{tab-examples2D.dat};
				\addlegendentry{$\|\widetilde u_{2,\varrho h}-\overline u_{2}\|_{L^2(Q)}$}
				\addplot table [col sep=&, y=err3, x=N]{tab-examples2D.dat};
				\addlegendentry{$\|\widetilde u_{3,\varrho h}-\overline u_{3}\|_{L^2(Q)}$}
				\addplot[
				domain = 1000:1000000,
				samples = 10,
				dashed,
				thin,
				blue,
				] {5*x^(-1)};
				\addlegendentry{$h^2\sim N^{-1}$}
				\addplot[
				domain = 10000:1000000,
				samples = 10,
				dashed,
				thin,
				red,
				] {1.5*x^(-1/4)};
				\addlegendentry{$h^{0.5}\sim N^{-1/3}$}
				\addplot[
				domain = 10000:1000000,
				samples = 10,
				dashed,
				thin,
				brown,
				] {5*x^(-3/4)};
				\addlegendentry{$h^{1.5}\sim N^{-3/4}$}
				\addplot[
				domain = 1000:10000,
				samples = 2,
				dashed,
				thin,
				brown,
				] {10^(-5)};
			\end{axis}
		\end{tikzpicture}
		\caption{Energy regularization, $\varrho =h^2$}
	\end{subfigure}
	\hfill
	\begin{subfigure}[b]{0.45\textwidth}
		\centering
		\begin{tikzpicture}[scale=0.7]
			\begin{axis}[
				xmode = log,
				ymode = log,
				xlabel={Number $N$ of space-time
					finite elements},
				ylabel=$\|\widetilde u_{i,\varrho h}-\overline u_{i}\|_{L^2(Q)}$,
				legend style={font=\tiny}, legend pos = south west]
				\addplot table [col sep=&, y=err1, x=N]{tab_examples_l2reg2D.dat};
				\addlegendentry{$\|\widetilde u_{1,\varrho h}-\overline u_{1}\|_{L^2(Q)}$}
				\addplot table [col sep=&, y=err2, x=N]{tab_examples_l2reg2D.dat};
				\addlegendentry{$\|\widetilde u_{2,\varrho h}-\overline u_{2}\|_{L^2(Q)}$}
				\addplot table [col sep=&, y=err3, x=N]{tab_examples_l2reg2D.dat};
				\addlegendentry{$\|\widetilde u_{3,\varrho h}-\overline u_{3}\|_{L^2(Q)}$}
				\addplot[
				domain = 10000:1000000,
				samples = 10,
				dashed,
				thin,
				blue,
				] {0.5*x^(-1)};
				\addlegendentry{$h^2\sim N^{-1}$}
				\addplot[
				domain = 10000:1000000,
				samples = 10,
				dashed,
				thin,
				red,
				] {1.5*x^(-1/4)};
				\addlegendentry{$h^{0.5}\sim N^{-1/3}$}
				\addplot[
				domain = 10000:1000000,
				samples = 10,
				dashed,
				thin,
				brown,
				] {5*x^(-3/4)};
				\addlegendentry{$h^{1.5}\sim N^{-3/4}$}
				\addplot[
				domain = 1000:10000,
				samples = 2,
				dashed,
				thin,
				brown,
				] {10^(-5)};
			\end{axis}
		\end{tikzpicture}
		\caption{$L^2$ regularization, $\varrho=h^4$}
	\end{subfigure}
	\caption{Convergence plots for the three different target functions
		$\overline u_i$, $i=1,2,3$ for the energy and the $L^2$
		regularization. }
	\label{fig:examples}
\end{figure}
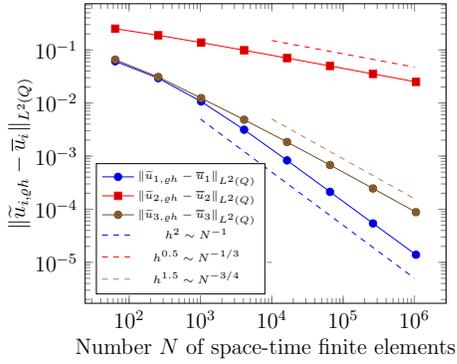
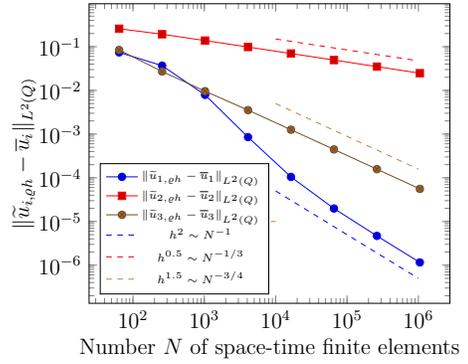

In order to check our theoretical findings, we consider three test examples
of different regularity for the target function $\overline u$, in the
space-time domain $Q:=(0,1)\times (0,1)\subset \mathbb{R}^{2}$. 
First we consider a smooth function
$\overline u_1\in C^2(\overline{Q})\cap H_{0;0,}^{1,1}(Q)$
given as 
\begin{align}\label{eq:smooth}
	\overline u_1(x,t) = \begin{cases}
		\frac 12(6t-3x-2)^3(3x-6t)^3,&\, x\leq t\text{ and }t-x\leq 2,\\
		0,&\text{else}. 
	\end{cases}
\end{align}
As a second target function we have the
piecewiese constant function $\overline u_2\in H^{1/2-\varepsilon}(Q)$,
$\varepsilon >0$, given as
\begin{align}\label{eq:discontinuous}
	\overline u_2(x,t) = \begin{cases} 1,& (x,t)\in (0.25,0.75)^2 \subset Q,\\
		0,& \text{else.}
	\end{cases}
\end{align}
Finally, we consider a piecewise bilinear
function $\overline u_3\in H_0^{3/2-\varepsilon}(Q)$,
$\varepsilon > 0$, defined as 
\begin{align}\label{eq:linear}
	\overline u_3(x,t) = \phi(x)\phi(t),\quad \phi(s)=\begin{cases}
		1, & s = 0.5,\\
		0, & s \not\in [0.25,0.75],\\
		\text{linear},& \text{else.}  
	\end{cases}
\end{align}
The numerical results for the energy regularization
\eqref{eq:matrix-form-energ-reg} with the optimal regularization
parameter $\varrho=h^2$, and for the $L^2$ regularization
\eqref{eq:matrix-form-l2reg} with $\varrho=h^4$,
are depicted in Fig.~\ref{fig:examples}, where we
observe optimal orders of convergence for each of the three examples,
as predicted by the theory. 

\begin{figure}[htbp]
	\centering
	\begin{subfigure}[b]{0.3\textwidth}
		\centering 
		\includegraphics[width=\textwidth]{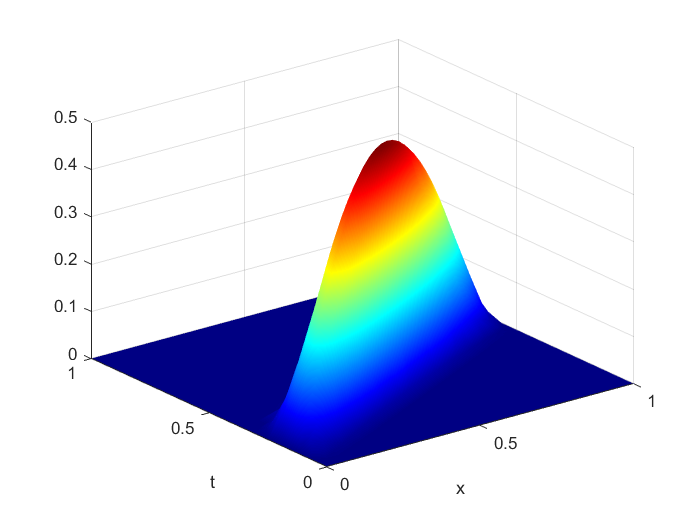}
		\caption{$\overline u_{1}$}
	\end{subfigure}
	\begin{subfigure}[b]{0.3\textwidth}
		\centering 
		\includegraphics[width=\textwidth]{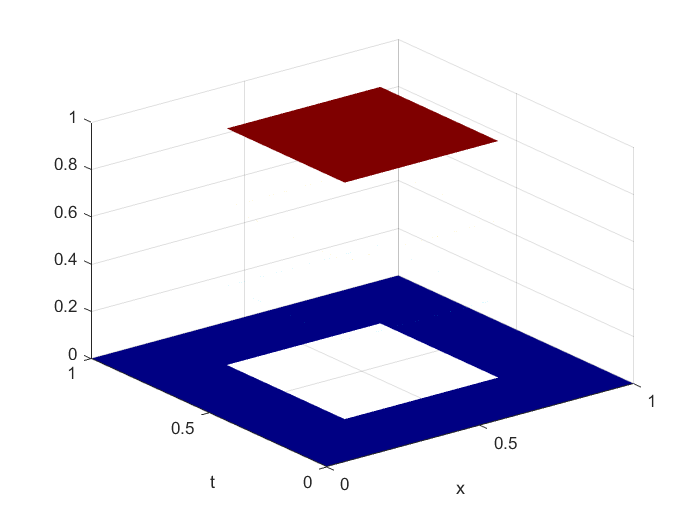}
		\caption{$\overline u_{2}$}
	\end{subfigure}
	\begin{subfigure}[b]{0.3\textwidth}
		\centering 
		\includegraphics[width=\textwidth]{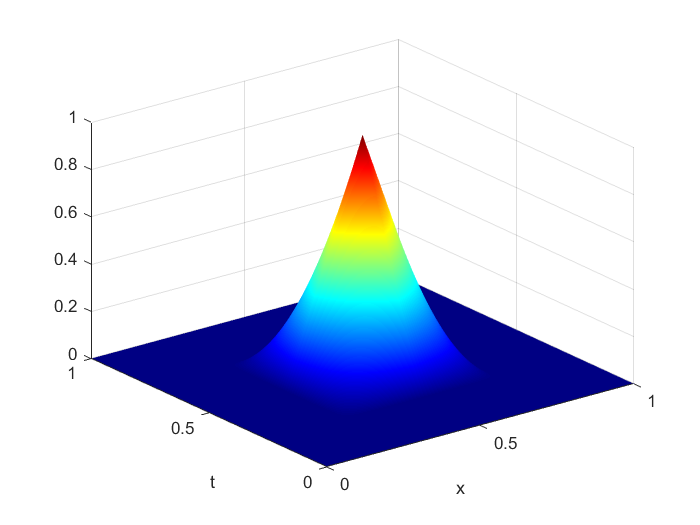}
		\caption{$\overline u_{3}$}
	\end{subfigure}
	\caption{Target functions $\overline u_{i}$, $i=1,2,3$.}
	\label{fig:adaptive-target-functions}
\end{figure}

%% 2D EXAMPLES END-------------------------------------%%%%%%%%%%%

\begin{figure}[htbp]
	\centering
	\begin{subfigure}[b]{0.3\textwidth}
		\centering 
		\includegraphics[width=\textwidth]{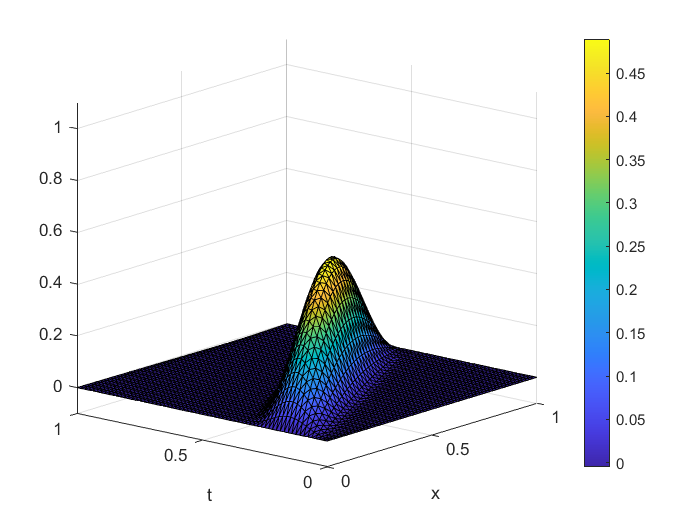}
	\end{subfigure}
	\hfill
	\begin{subfigure}[b]{0.3\textwidth}
		\centering 
		\includegraphics[width=\textwidth]{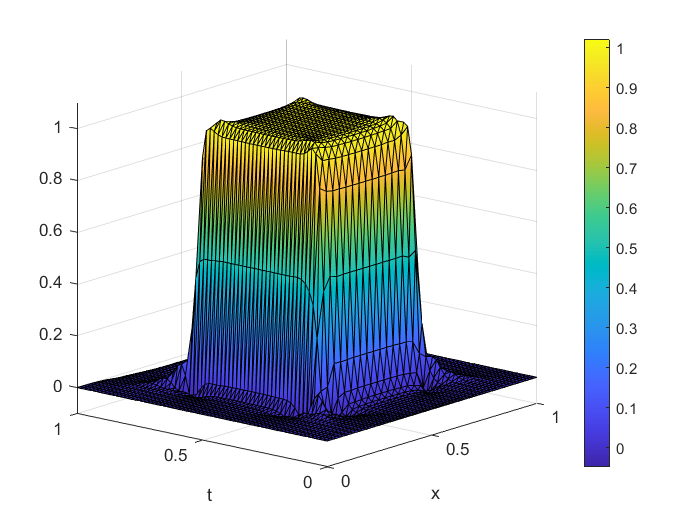}
	\end{subfigure}
	\hfill      
	\begin{subfigure}[b]{0.3\textwidth}
		\centering 
		\includegraphics[width=\textwidth]{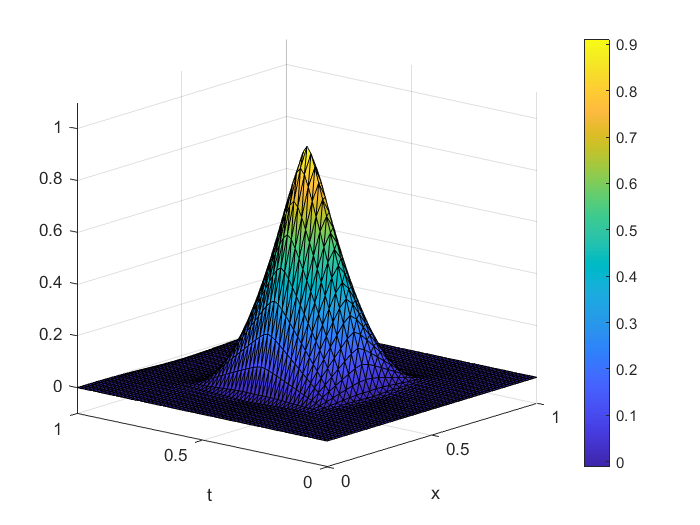}
	\end{subfigure}
	\\
	Computed states in the case of energy regularization
	\\
	\begin{subfigure}[b]{0.3\textwidth}
		\centering 
		\includegraphics[width=\textwidth]{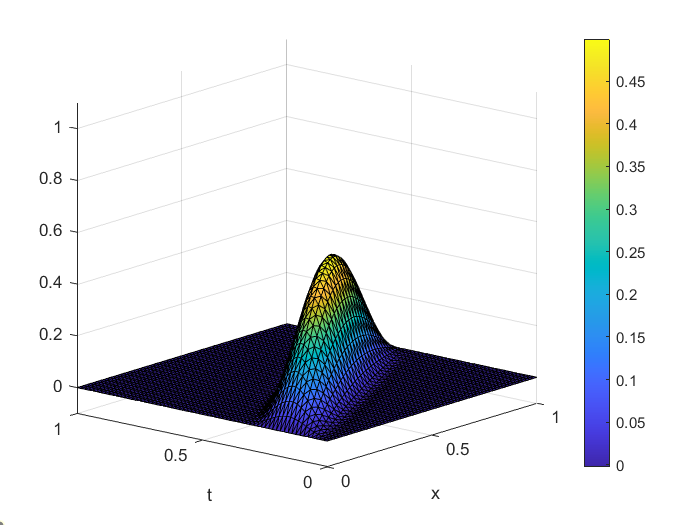}
	\end{subfigure}
	\hfill
	\begin{subfigure}[b]{0.3\textwidth}
		\centering 
		\includegraphics[width=\textwidth]{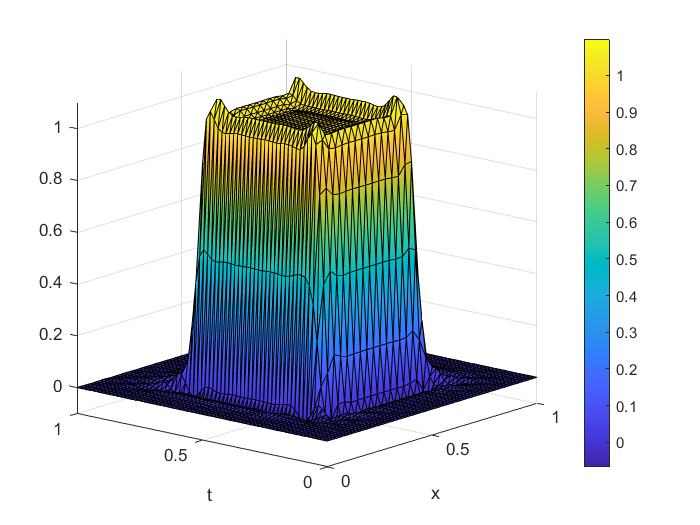}
	\end{subfigure}
	\hfill
	\begin{subfigure}[b]{0.3\textwidth}
		\centering 
		\includegraphics[width=\textwidth]{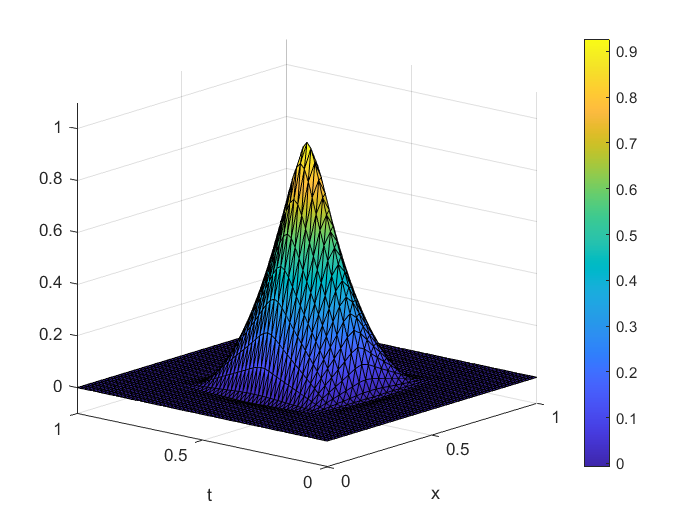}
	\end{subfigure}
	\\
	Computed states in the case of $L^2$ regularization
	\\
	\begin{subfigure}[b]{0.3\textwidth}
		\centering 
		\includegraphics[width=\textwidth]{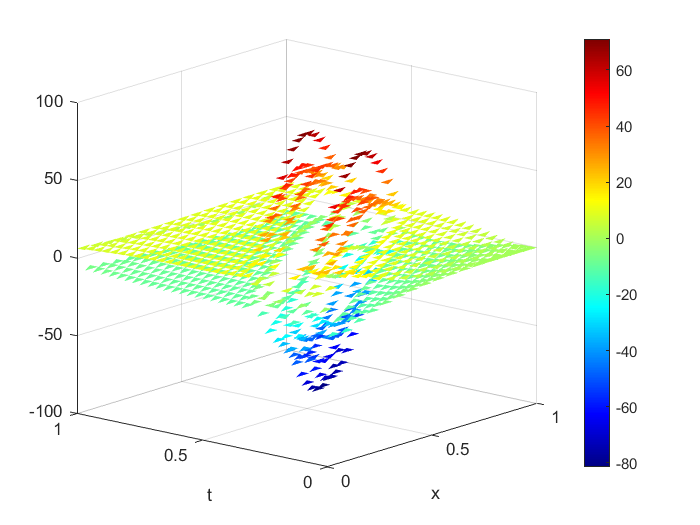}
	\end{subfigure}
	\hfill
	\begin{subfigure}[b]{0.3\textwidth}
		\centering 
		\includegraphics[width=\textwidth]{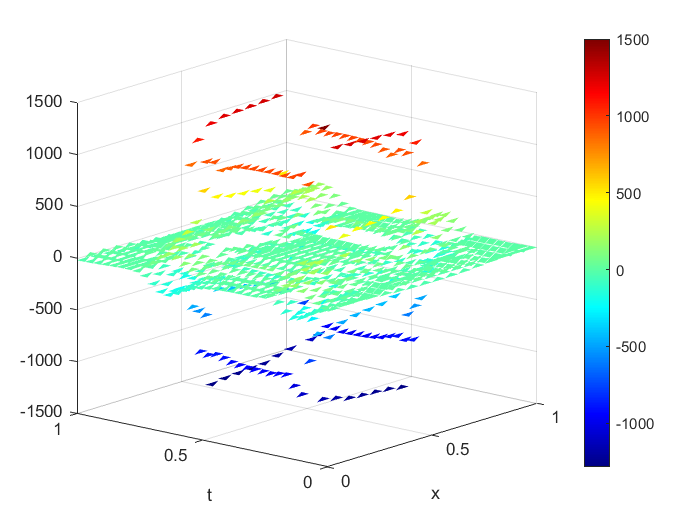}
	\end{subfigure}
	\hfill
	\begin{subfigure}[b]{0.3\textwidth}
		\centering 
		\includegraphics[width=\textwidth]{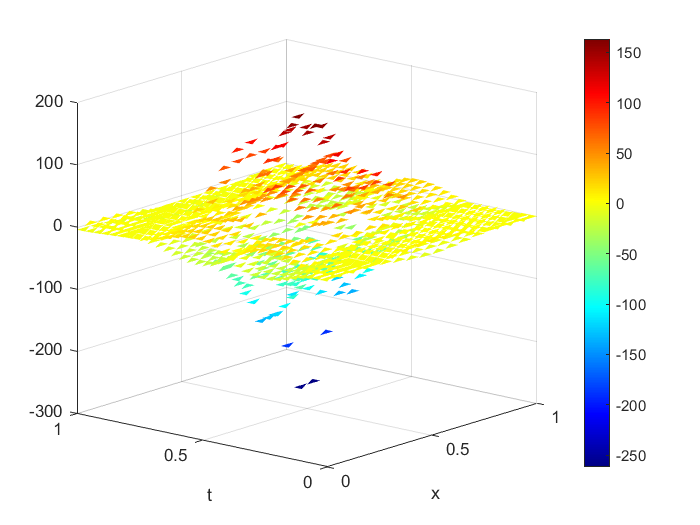}
	\end{subfigure}
	\\
	Reconstructed controls in the case of energy regularization
	\\
	\begin{subfigure}[b]{0.3\textwidth}
		\centering 
		\includegraphics[width=\textwidth]{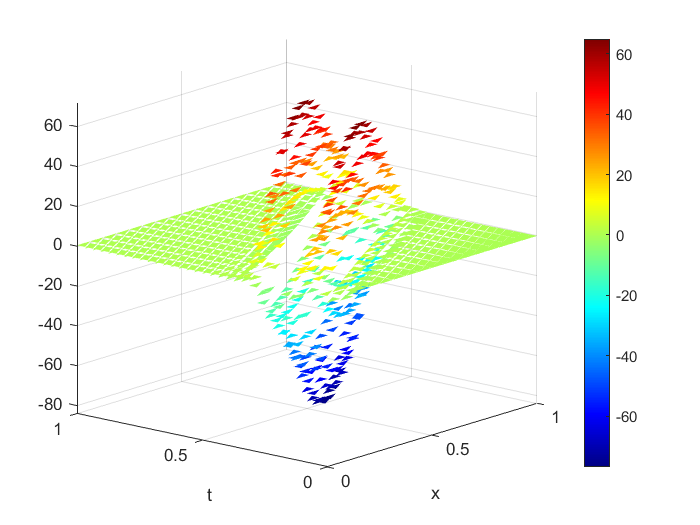}
	\end{subfigure}
	\hfill
	\begin{subfigure}[b]{0.3\textwidth}
		\centering 
		\includegraphics[width=\textwidth]{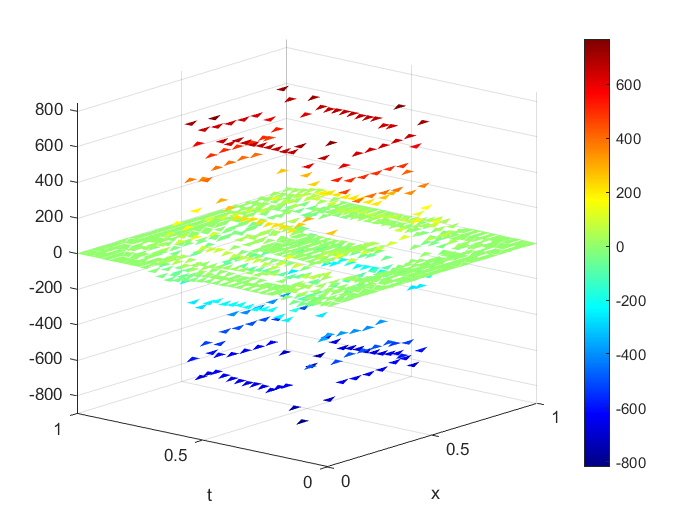}
	\end{subfigure}
	\hfill
	\begin{subfigure}[b]{0.3\textwidth}
		\centering 
		\includegraphics[width=\textwidth]{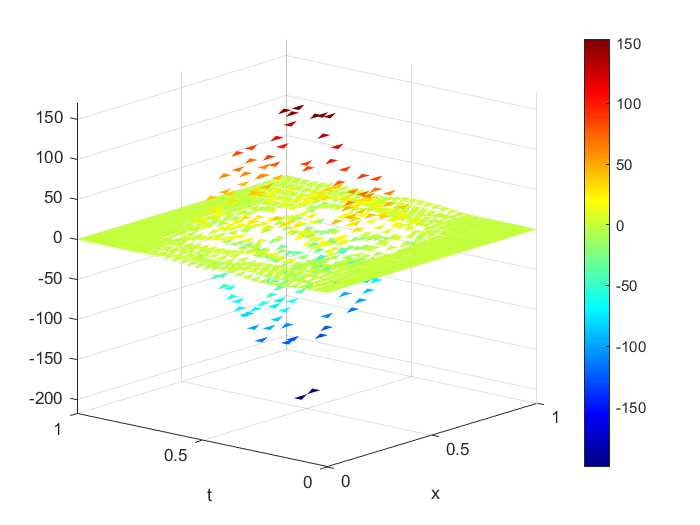}
	\end{subfigure}
	\\
	Reconstructed controls in the case of $L^2$ regularization
	\caption{Comparison of the computed states and the reconstructed
		controls for energy and $L^2$ regularizations on level 3, with
		4096 elements and 1984 degrees of freedom.}
	\label{fig:comparison-l2-energ-reg}
\end{figure}

\noindent
Although the convergence rates and errors for both approaches, the energy
and the $L^2$ regularization method, seem to be comparable, we observe a
difference in the behaviour of the discontinuous solution
$\widetilde u_{2,\varrho h}$, see Fig.~\ref{fig:comparison-l2-energ-reg}.
This is due to the additional regularity $z_\varrho\in H_{0;,0}^{1,1}(Q)$
which we gain when considering the control in $Z=L^2(Q)$.

\begin{remark}
	The choice of $A$ enforces homogeneous Neumann
	conditions at the origin $t=0$, while we have homogeneous
	Dirichlet boundary conditions elsewhere. Due to this change
	in the boundary conditions we may have a reduced regularity
	for the solution of the space-time Poisson equation, see
	Corollary \ref{corollary regularization}. This results in a
	reduced order of convergence, as observed for the target
	$\overline u_4(x,t)=t\sin(\pi t)\sin(\pi x)$, for $(x,t)\in (0,1)^2$,
	where the solution of the energy regularization
	\eqref{eq:matrix-form-energ-reg} converges with a rate
	$1.5$ instead of $2$, see
	Table \ref{tab:slower-convergence}.
	To regain optimal rates, there are three possible remedies. First,
	one might choose $\varrho =h^3$ for the energy regularization. Then
	the term is penalized strong enough to ensure optimal orders of
	convergence. Second, one might compute the solution on an enlarged
	domain, embedding the target function such that in a neighborhood of
	$t=0$ and $t=T$ the function is constant zero, as we have done in
	our examples $\overline{u}_i$, $i=1,2,3$.
	Then $p_\varrho$ will (approximately) fulfil
	the homogeneous Neumann condition. A third possibility is to adaptively
	refine the mesh and resolve the singularities. This will be discussed
	in the next section. Note, that for the $L^2$ regularization approach
	this effect does not occur, since the operator
	$A=\text{id}:L^2(Q)\to L^2(Q)$ does not enforce any initial condition.
\end{remark}  

\begin{table}[htbp]
	\centering
	\begin{filecontents*}{tab_smooth_less-convergence2D.dat}
level&M&N&h_min&h_max&rho&err1&eoc1
0&24&64&0.125&0.125&0.015625&0.024620526784637&0
1&112&256&0.0625&0.0625&0.00390625&0.00774614213528852&1.66831175130093
2&480&1024&0.03125&0.03125&0.0009765625&0.00269562925497814&1.52285588565067
3&1984&4096&0.015625&0.015625&0.000244140625&0.00096343345604377&1.48436515950068
4&8064&16384&0.0078125&0.0078125&6.103515625e-05&0.00034448172342557&1.48375758123344
5&32512&65536&0.00390625&0.00390625&1.52587890625e-05&0.000122676704236635&1.48956612876573
6&130560&262144&0.001953125&0.001953125&3.814697265625e-06&4.35480872157422e-05&1.49418005648654
7&523264&1048576&0.0009765625&0.0009765625&9.5367431640625e-07&1.54293051239311e-05&1.49693626202343
\end{filecontents*}
\pgfplotstabletypeset[col sep=&,
	columns={level,M,N,h_min,rho,err1,eoc1},
	columns/level/.style={column name=Level},
	columns/M/.style={fixed, column name=DoFs},
	columns/N/.style={fixed, column name=$N$},
	columns/h_min/.style={fixed, fixed zerofill, precision=3,
		column name=$h$},
	columns/rho/.style={sci, sci zerofill, precision=2, column
		name=$\varrho\, (h^2)$},
	columns/err1/.style={sci,sci zerofill,column name=$\|\widetilde u_{4,\varrho h}-\overline u_4\|_{L^2(Q)}$},
	columns/eoc1/.style={fixed, fixed zerofill, precision=2, column name=eoc},
	every head row/.style={before row=\toprule , after row=\midrule},
	every last row/.style={after row=\bottomrule}]
	{tab_smooth_less-convergence2D.dat}
	\caption{Errors and orders of convergence for
		$\overline u_4(x,t)=t\sin(t\pi)\sin(x\pi)$ in the case of an
		uniform refinement strategy with $\varrho= h^2$.}
	\label{tab:slower-convergence}
\end{table}  

\noindent
Since in optimal control theory we are mainly interested in
the control $z$, rather than in the computed state
$\widetilde{u}_{\varrho h}$, we are going to reconstruct the control
$z_\varrho$ in a post-processing step
when solving \eqref{eq:VF_control-discrete}. We introduce the
finite element space
$Z_H:=S_H^0(Q)=\text{span}\{\phi_r\}_{r=1}^{N_H}\subset [H_{0;,0}^{1,1}(Q)]^*$
of piecewise constant basis functions $\phi_r$. When using
$Y_h:=S_h^1(Q)\cap Y$, \eqref{eq:VF_control-discrete} is equivalent to
the linear system of algebraic equations
to find $\underline \psi\in\mathbb{R}^{M_Y}\leftrightarrow\psi_h\in Y_h$ and
$\underline z\in\mathbb{R}^{N_H} \leftrightarrow \widetilde z_{\varrho H}\in Z_H$
such that
\begin{align*}
	\begin{pmatrix}
		A_h & P_{hH}^\top\\ P_{hH} & 0
	\end{pmatrix}\begin{pmatrix}
		\underline \psi\\\underline z
	\end{pmatrix}=\begin{pmatrix}
		B_h \underline u\\0
	\end{pmatrix},
\end{align*}  
with matrices $B_h$ as above and $$ P_{hH}[r,j]=\langle \psi_j,\phi_r  \rangle_{L^2(Q)},\quad r=1,\ldots,N_H, \, j = 1,\ldots, M_Y. $$	 
Resolving the system for $\underline z$ gives 
\begin{equation}\label{eq:control-reconstruction-matrix}
	\underline z = (P_{hH} A_h^{-1}P_{hH}^\top)^{-1}P_{hH}A_h^{-1}B_h\underline u.
\end{equation}
In Fig.~\ref{fig:comparison-l2-energ-reg} we also present the reconstructed
controls for both the energy and the $L^2$ regularization approach.

\subsection{Adaptive refinement}
\begin{figure}[h]
	\centering
	\begin{subfigure}[b]{0.3\textwidth}
		\centering 
		\includegraphics[width=\textwidth]{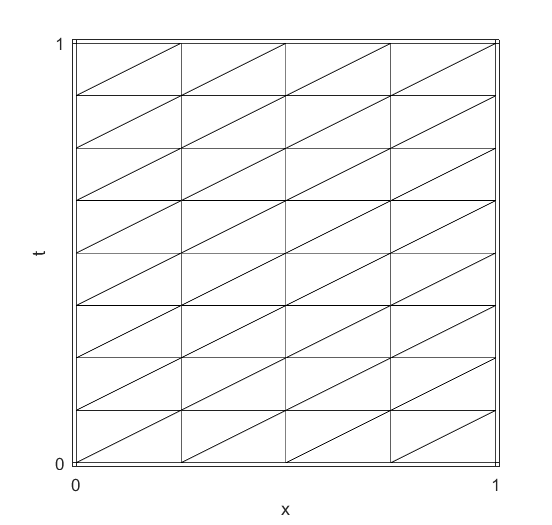}
		\caption{Level 0}
	\end{subfigure}
	\hfill
	\begin{subfigure}[b]{0.3\textwidth}
		\centering 
		\includegraphics[width=\textwidth]{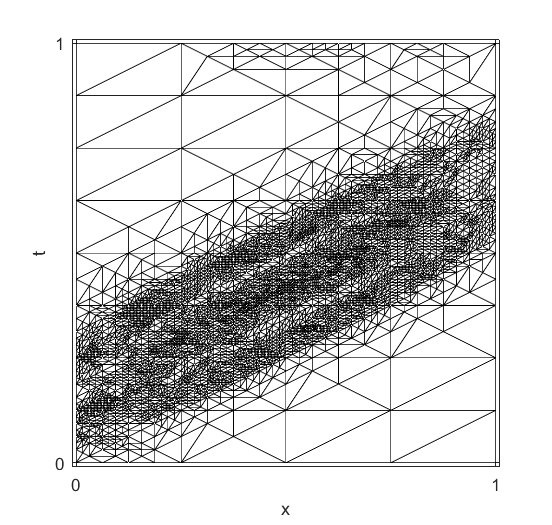}
		\caption{$\overline u_1$, Level 10}
	\end{subfigure}
	\hfill
	\begin{subfigure}[b]{0.3\textwidth}
		\centering 
		\includegraphics[width=\textwidth]{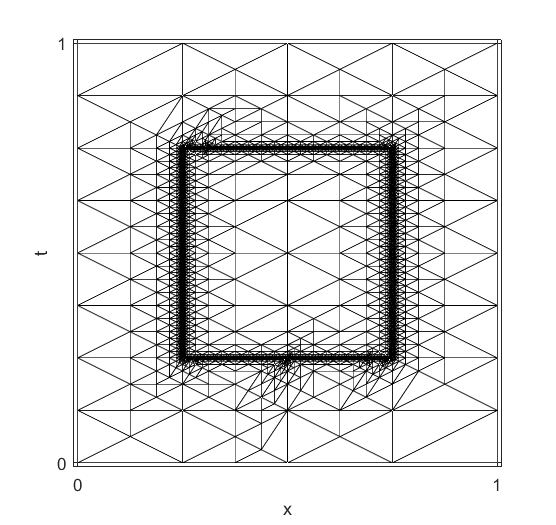}
		\caption{$\overline u_2$, Level 7}
	\end{subfigure}
	\caption{Initial mesh and adaptively refined meshes for the target functions $\overline u_{i}$, $i=1,2$.}
	\label{fig:initial-and-refined-meshes}
\end{figure}

In this section we present some examples for an adaptive space-time
refinement strategy for the energy regularization
\eqref{eq:matrix-form-energ-reg}. We will apply an adaptive refinement
strategy using D\"orfler marking \cite{Doerfler:1996} with the refinement
indicator
$ \eta_\ell=\|\widetilde u_{\varrho h}-\overline u\|_{L^2(\tau_\ell)}$
on each simplicial space-time finite element $\tau_\ell$, $\ell=1,\ldots,N$.
With this choice we see that the approximation error fufills
\[
\|\widetilde u_{\varrho h}-\overline u\|_{L^2(Q)}^2 =
\sum_{\ell=1}^N \eta_\ell^2.
\]
We will refine all elements $\tau_k$ that satisfy
\[
\eta_k\geq \theta \max_{\ell=1,\ldots,N}\eta_\ell,
\]
with $\theta = 0.5$. The initial mesh with $64$ elements and $24$ degrees
of freedom (DoFs) and the resulting adaptively refined meshes for the
target functions $\overline u_{1}$ at level $10$ with $8824$ elements and
$4389$ DoFs and for $\overline u_{2}$ at level $7$ with $15159$ elements
and $7571$ DoFs are shown in Fig.~\ref{fig:initial-and-refined-meshes}.
In Tables \ref{tab:comparison-uniform-adaptive-u1} and
\ref{tab:comparison-uniform-adaptive-u2} we present a comparison of the
errors of the adaptive refinement strategy against the errors of the
uniform refinement at levels with comparably many elements for both
target functions. We clearly see, that in both cases considerably less
elements are needed to achieve errors of the same order.

\begin{table}[h]
	\begin{tabular}{|l|rlr|l|rlr|}
		\hline
		\multirow{2}{*}{L}&\multicolumn{3}{c|}{Adaptive} & \multirow{2}{*}{L}&
		\multicolumn{3}{c|}{Uniform} \\  \cline{2-4} \cline{6-8}
		&\#DoFs&$\|\widetilde{u}_{1,\varrho h}-\overline{u}_1\|_{L^2(Q)}$ & $\varrho=h_{min}^2$&&\#DoFs & $\|\widetilde{u}_{1,\varrho h}-\overline{u}_1\|_{L^2(Q)}$ &$\varrho=h^2$\\
		\hline
		0&$24$ &$6.12415\cdot10^{-2}$&$1.56\cdot10^{-2}$ &						0&$24$		&$6.12415\cdot10^{-2}$ &$1.56\cdot10^{-2}$ \\
		2&$101$&$1.26895\cdot10^{-2}$&$9.77\cdot10^{-4}$ &						1&$112$		&$2.94242\cdot10^{-2}$ &$3.91\cdot10^{-3}$ \\
		5&$399$&$2.05047\cdot10^{-3}$&$1.22\cdot10^{-4}$ &						2&$480$		&$1.06888\cdot10^{-2}$ &$9.77\cdot10^{-4}$ \\
		8&$1654$&$5.13791\cdot10^{-4}$&$3.05\cdot10^{-5}$ &						3&$1984$	&$3.14290\cdot10^{-3}$ &$2.44\cdot10^{-4}$ \\
		13&$8266$&$6.61367\cdot10^{-5}$&$1.91\cdot10^{-6}$ &					4&$8064$	&$8.31332\cdot10^{-4}$ &$6.10\cdot10^{-5}$ \\
		18&$39821$&$1.19888\cdot10^{-5}$&$1.19\cdot10^{-7}$ &					5&$32512$	&$2.12588\cdot10^{-4}$ &$1.53\cdot10^{-5}$ \\
		22&$162774$&$2.86782\cdot10^{-6}$&$2.98\cdot10^{-8}$ &					6&$130560$	&$5.41480\cdot10^{-5}$ &$3.82\cdot10^{-6}$ \\
		25&$377896$&$1.07848\cdot10^{-6}$&$7.45\cdot10^{-9}$ &					7&$523264$	&$1.39319\cdot10^{-5}$ &$9.54\cdot10^{-7}$  \\
		26&$636878$&$7.03678\cdot10^{-7}$&$1.86\cdot10^{-9}$ &		& 			&						& \\
		\hline
	\end{tabular}
	\caption{Comparison of the uniform refinement to the adaptive refinement strategy for levels (L) with comparably many elements  for the target function $\overline u_{1}$ for the energy regularization \eqref{eq:matrix-form-energ-reg}. } 
	\label{tab:comparison-uniform-adaptive-u1}
\end{table}

\begin{table}[h]
	\begin{tabular}{|l|rlr|l|rlr|}
		\hline
		\multirow{2}{*}{L}&\multicolumn{3}{c|}{Adaptive} & \multirow{2}{*}{L}&
		\multicolumn{3}{c|}{Uniform} \\  \cline{2-4} \cline{6-8}
		&\#DoFs&$\|\widetilde{u}_{2,\varrho h}-\overline{u}_2\|_{L^2(Q)}$ & $\varrho=h_{min}^2$&&\#DoFs & $\|\widetilde{u}_{2,\varrho h}-\overline{u}_2\|_{L^2(Q)}$ &$\varrho=h^2$\\
		\hline
		0&$24$&$2.50691\cdot10^{-1}$ &$1.56\cdot10^{-2}$ &						0&$24$		&$2.50691\cdot10^{-1}$ &$1.56\cdot10^{-2}$ \\
		2&$198$&$1.36350\cdot10^{-1}$&$9.77\cdot10^{-4}$ &						1&$112$		&$1.88590\cdot10^{-1}$ &$3.91\cdot10^{-3}$ \\
		3&$435$&$9.74050\cdot10^{-2}$&$2.44\cdot10^{-4}$ &						2&$480$		&$1.37373\cdot10^{-1}$ &$9.77\cdot10^{-4}$ \\
		5&$1895$&$4.92039\cdot10^{-2}$&$1.53\cdot10^{-5}$ &					3&$1984$	&$9.85712\cdot10^{-2}$ &$2.44\cdot10^{-4}$ \\
		7&$7571$&$2.46665\cdot10^{-2}$&$9.54\cdot10^{-7}$ &					4&$8064$	&$7.02300\cdot10^{-2}$ &$6.10\cdot10^{-5}$ \\
		9&$30027$&$1.23436\cdot10^{-2}$&$5.96\cdot10^{-8}$ &					5&$32512$	&$4.98503\cdot10^{-2}$ &$1.53\cdot10^{-5}$ \\
		11&$119554$&$6.17867\cdot10^{-3}$&$3.73\cdot10^{-9}$ &					6&$130560$	&$3.53171\cdot10^{-2}$ &$3.82\cdot10^{-6}$ \\
		13&$477542$&$3.09069\cdot10^{-3}$&$2.3\cdot10^{-10}$ &				7&$523264$	&$2.49969\cdot10^{-2}$ &$9.54\cdot10^{-7}$  \\
		14&$957389$&$2.18324\cdot10^{-3}$&$5.8\cdot10^{-11}$ &		         & 			&						& \\
		\hline
	\end{tabular}
	\caption{Comparison of the uniform refinement to the adaptive refinement strategy for levels (L) with comparably many elements  for the target function $\overline u_{2}$ for the energy regularization \eqref{eq:matrix-form-energ-reg}. } 
	\label{tab:comparison-uniform-adaptive-u2}
\end{table}

\noindent
All comuptations were carried out with Matlab using a sparse direct
solver. For the adaptive refinement strategy the package
from \cite{FunkenPraetoriusWissgott:2011} was adapted suitably. 

\section{Conclusions and outlook}
We have introduced and investigated a space-time finite element method for
distributed optimal control problems for the wave equation with energy
regularization. In particular, we have shown $L^2(Q)$ error estimates
between the desired state $\overline{u}$ and the computable discrete
solution $\widetilde{u}_{\varrho h}$, with respect to the regularity of
the target function. It has been proven that in this case the choice
$\varrho=h^2$ delivers optimal orders of convergence, and the findings have
been supported by several numerical examples. Moreover, we compared the
results to the case where a $L^2(Q)$ regularization is used. Furthermore,
we proposed an adaptive finite element strategy and
presented its performance for target functions with different regularities,
where we observed that considerably less elements are needed for a
comparable error than in the case of an uniform refinement. 

The system matrices, for both the $L^2$ and the energy
regularization approach, are positive definite but skew-symmetric, or
alternatively, symmetric but indefinite. Thus it is of highest interest to
develop robust iterative solvers as already done in the elliptic and
parabolic case
\cite{LangerLoescherSteinbachYang:2022a,UL:LangerSteinbachYang:2022a,LLSY2:LangerSteinbachYang:2022b}.
This will then also allow an efficient solution of related optimal
control problems in two and three space dimensions.
Moreover,
for discontinuous targets and an adaptive finite element scheme, it would be
sensible to consider a relaxation parameter $\varrho=\varrho(x,t)$ that is
locally varying, to have a better resolution of the control defined on the
adaptive mesh.
This has already been done for the optimal control problem subject
to the Poisson equation \cite{LLSY:2022variable}. 
In addition, in order to be of practical interest,
the consideration of control and/or state constraints can be considered
within the abstract framework as done in \cite{GanglLoescherSteinbach:2022}.

\bigskip

\noindent
{\bf Acknowledgment:} The authors would like to thank U.~Langer and
F.~Tr\"oltzsch for the fruitful discussions and their helpful
comments during their visit to TU Graz in October 2022.

\bibliography{wave}

\begin{thebibliography}{10}

\bibitem{Adams:1975}
R.~A. Adams.
\newblock {\em Sobolev Spaces}.
\newblock Academic Press, New York, London, 1975.

\bibitem{AlbaniCezaroZubelli:2016}
V.~Albani, A.~De~Cezaro, and J.~P. Zubelli.
\newblock On the choice of the {T}ikhonov regularization parameter and the
  discretization level: a discrepancy-based strategy.
\newblock {\em Inverse Probl. Imaging}, 10(1):1--25, 2016.

\bibitem{LSTY:BabuskaAziz:1972a}
I.~Babu\v{s}ka and A.~Aziz.
\newblock Survey lectures on the mathematical foundation of the finite element
  method.
\newblock In {\em The Mathematical Foundations of the Finite Element Method
  with Applications to Partial Differential Equations}, pages 1--359, New York,
  1972. Academic Press.

\bibitem{BangerthGeigerRannacher:2010}
W.~Bangerth, M.~Geiger, and R.~Rannacher.
\newblock Adaptive {G}alerkin finite element methods for the wave equation.
\newblock {\em Comput. Methods Appl. Math.}, 10(1):3--48, 2010.

\bibitem{BrennerScott:2008}
S.~C. Brenner and L.~R. Scott.
\newblock {\em The mathematical theory of finite element methods}, volume~15 of
  {\em Texts in Applied Mathematics}.
\newblock Springer, New York, third edition, 2008.

\bibitem{Doerfler:1996}
W.~D\"{o}rfler.
\newblock A convergent adaptive algorithm for {P}oisson's equation.
\newblock {\em SIAM J. Numer. Anal.}, 33(3):1106--1124, 1996.

\bibitem{DoerflerFindeisenWienersZiegler:2019}
W.~D\"{o}rfler, S.~Findeisen, C.~Wieners, and D.~Ziegler.
\newblock Parallel adaptive discontinuous {G}alerkin discretizations in space
  and time for linear elastic and acoustic waves.
\newblock In {\em Space-time methods---applications to partial differential
  equations}, volume~25 of {\em Radon Ser. Comput. Appl. Math.}, pages 61--88.
  de Gruyter, Berlin, 2019.

\bibitem{EnglHankeNeubauer:1996}
H.~W. Engl, M.~Hanke, and A.~Neubauer.
\newblock {\em Regularization of inverse problems}, volume 375 of {\em
  Mathematics and its Applications}.
\newblock Kluwer Academic Publishers Group, Dordrecht, 1996.

\bibitem{ErnestiWieners:2019}
J.~Ernesti and C.~Wieners.
\newblock A space-time discontinuous {P}etrov-{G}alerkin method for acoustic
  waves.
\newblock In {\em Space-time methods---applications to partial differential
  equations}, volume~25 of {\em Radon Ser. Comput. Appl. Math.}, pages 89--115.
  de Gruyter, Berlin, 2019.

\bibitem{FunkenPraetoriusWissgott:2011}
S.~Funken, D.~Praetorius, and P.~Wissgott.
\newblock Efficient implementation of adaptive {P}1-{FEM} in {M}atlab.
\newblock {\em Comput. Methods Appl. Math.}, 11(4):460--490, 2011.

\bibitem{GanglLoescherSteinbach:2022}
P.~Gangl, R.~L\"oscher, and O.~Steinbach.
\newblock Regularization and finite element error estimates for optimal control
  problems with energy regularization and state or control constraints.
\newblock In preparation, 2022.

\bibitem{GlowinskiKintonWheeler:1989}
R.~Glowinski, W.~Kinton, and M.~F. Wheeler.
\newblock A mixed finite element formulation for the boundary controllability
  of the wave equation.
\newblock {\em Internat. J. Numer. Methods Engrg.}, 27(3):623--635, 1989.

\bibitem{GugatKeimerLeugering:2009}
M.~Gugat, A.~Keimer, and G.~Leugering.
\newblock Optimal distributed control of the wave equation subject to state
  constraints.
\newblock {\em Z. Angew. Math. Mech.}, 89(6):420--444, 2009.

\bibitem{HulbertHughes:1990}
G.~M. Hulbert and T.~J.~R. Hughes.
\newblock Space-time finite element methods for second-order hyperbolic
  equations.
\newblock {\em Comput. Methods Appl. Mech. Engrg.}, 84(3):327--348, 1990.

\bibitem{Isakov:2017}
V.~Isakov.
\newblock {\em Inverse problems for partial differential equations}, volume 127
  of {\em Applied Mathematical Sciences}.
\newblock Springer, Cham, third edition, 2017.

\bibitem{KroenerKunischVexler:2011}
A.~Kr\"{o}ner, K.~Kunisch, and B.~Vexler.
\newblock Semismooth {N}ewton methods for optimal control of the wave equation
  with control constraints.
\newblock {\em SIAM J. Control Optim.}, 49(2):830--858, 2011.

\bibitem{Ladyzhenskaya:1985}
O.~A. Ladyzhenskaya.
\newblock {\em The boundary value problems of mathematical physics}, volume~49
  of {\em Applied Mathematical Sciences}.
\newblock Springer-Verlag, New York, 1985.

\bibitem{LLSY:2022variable}
U.~Langer, R.~L\"oscher, O.~Steinbach, and H.~Yang.
\newblock An adaptive finite element method for distributed elliptic optimal
  control problems with variable energy regularization.
\newblock arXiv:2209.08811, 2022.

\bibitem{LangerLoescherSteinbachYang:2022a}
U.~Langer, R.~L\"oscher, O.~Steinbach, and H.~Yang.
\newblock Robust finite element discretization and solvers for distributed
  elliptic optimal control problems.
\newblock arXiv:2207.04664, 2022.

\bibitem{LangerSteinbachTroeltzschYang:2021}
U.~Langer, O.~Steinbach, F.~Tr\"{o}ltzsch, and H.~Yang.
\newblock Space-time finite element discretization of parabolic optimal control
  problems with energy regularization.
\newblock {\em SIAM J. Numer. Anal.}, 59(2):675--695, 2021.

\bibitem{LangerSteinbachTroeltzschYang:2022}
U.~Langer, O.~Steinbach, F.~Tr\"{o}ltzsch, and H.~Yang.
\newblock Unstructured space-time finite element methods for optimal sparse
  control of parabolic equations.
\newblock In {\em Optimization and control for partial differential
  equations---uncertainty quantification, open and closed-loop control, and
  shape optimization}, volume~29 of {\em Radon Ser. Comput. Appl. Math.}, pages
  167--188. de Gruyter, Berlin, 2022.

\bibitem{UL:LangerSteinbachYang:2022a}
U.~Langer, O.~Steinbach, and H.~Yang.
\newblock Robust space-time finite element error estimates for parabolic
  distributed optimal control problems with energy regularization.
\newblock arXiv:2206.06455, 2022.

\bibitem{LLSY2:LangerSteinbachYang:2022b}
U.~Langer, O.~Steinbach, and H.~Yang.
\newblock Robust discretization and solvers for elliptic optimal control
  problems with energy regularization.
\newblock {\em Comput. Methods Appl. Math.}, 22:97--111, 2022.

\bibitem{LionsVol1:1971}
J.-L. Lions.
\newblock {\em Optimal control of systems governed by partial differential
  equations}.
\newblock Die Grundlehren der mathematischen Wissenschaften, Band 170.
  Springer-Verlag, New York-Berlin, 1971.

\bibitem{LionsMagenesVol1:1968}
J.-L. Lions and E.~Magenes.
\newblock {\em Probl\`emes aux limites non homog\`enes et applications, Vol.
  1}, volume~17 of {\em Travaux et Recherches Math\'{e}matiques}.
\newblock Dunod, Paris, 1968.

\bibitem{McLean:2000}
W.~McLean.
\newblock {\em Strongly Elliptic Systems and Boundary Integral Equations}.
\newblock Cambridge University Press, Cambridge, UK, 2000.

\bibitem{MeidnerVexler:2007}
D.~Meidner and B.~Vexler.
\newblock Adaptive space-time finite element methods for parabolic optimization
  problems.
\newblock {\em SIAM J. Control Optim.}, 46(1):116--142, 2007.

\bibitem{MontanerMunch:2019}
S.~Montaner and A.~M\"{u}nch.
\newblock Approximation of controls for linear wave equations: a first order
  mixed formulation.
\newblock {\em Math. Control Relat. Fields}, 9(4):729--758, 2019.

\bibitem{LLSY2:NeumuellerSteinbach:2021a}
M.~Neum\"uller and O.~Steinbach.
\newblock Regularization error estimates for distributed control problems in
  energy spaces.
\newblock {\em Math. Methods Appl. Sci.}, 44:4176--4191, 2021.

\bibitem{LSTY:Necas:1962a}
J.~Ne\v{c}as.
\newblock Sur une m\'{e}thode pour r\'{e}soudre les \'{e}quations aux
  d\'{e}riv\'{e}es partielles du type elliptique, voisine de la variationnelle.
\newblock {\em Ann. Scuola Norm. Sup. Pisa}, 16(4):305--326, 1962.

\bibitem{KunischPeralta:2022}
G.~Peralta and K.~Kunisch.
\newblock Mixed and hybrid {P}etrov-{G}alerkin finite element discretization
  for optimal control of the wave equation.
\newblock {\em Numer. Math.}, 150(2):591--627, 2022.

\bibitem{Poerner:2018}
F.~P\"orner.
\newblock {\em Regularization Methods for Ill-Posed Optimal Control Problems}.
\newblock PhD thesis, Julius--Maximilians--Universit\"at W\"urzburg, 2018.

\bibitem{SteinbachZank:2020}
O.~Steinbach and M.~Zank.
\newblock Coercive space-time finite element methods for initial boundary value
  problems.
\newblock {\em Electron. Trans. Numer. Anal.}, 52:154--194, 2020.

\bibitem{SteinbachZank:2022}
O.~Steinbach and M.~Zank.
\newblock A generalized inf-sup stable variational formulation for the wave
  equation.
\newblock {\em J. Math. Anal. Appl.}, 505(1), 2022.
\newblock Paper No. 125457, 24 pp.

\bibitem{Troeltzsch:2010}
F.~Tr\"{o}ltzsch.
\newblock {\em Optimal control of partial differential equations}, volume 112
  of {\em Graduate Studies in Mathematics}.
\newblock American Mathematical Society, Providence, RI, 2010.

\bibitem{Zuazua:2005}
E.~Zuazua.
\newblock Propagation, observation, and control of waves approximated by finite
  difference methods.
\newblock {\em SIAM Rev.}, 47(2):197--243, 2005.

\end{thebibliography}
\bibliographystyle{abbrv}

\end{document}